\documentclass[11pt,a4paper,reqno]{amsart}
\usepackage[margin=28mm]{geometry}
\usepackage[T1]{fontenc}
\usepackage[utf8]{inputenc}
\usepackage{lmodern}
\usepackage{amsmath,amssymb,mathtools}
\usepackage{mathrsfs}
\usepackage{microtype}
\usepackage{xcolor}
\definecolor{linkblue}{RGB}{24,57,91}
\usepackage[colorlinks=true,linkcolor=linkblue,citecolor=linkblue,urlcolor=linkblue]{hyperref}
\usepackage{bookmark}
\usepackage{enumitem}
\usepackage{comment}

\numberwithin{equation}{section}
\newtheorem{theorem}{Theorem}[section]
\newtheorem{proposition}[theorem]{Proposition}
\newtheorem{lemma}[theorem]{Lemma}

\theoremstyle{definition}

\theoremstyle{remark}
\newtheorem{remark}[theorem]{Remark}
\theoremstyle{plain}
\newtheorem{mainresult}{Theorem}

\newtheorem{maincorollary}[mainresult]{Corollary}
\newcommand{\R}{\mathbb R}
\newcommand{\RP}{\mathbb{RP}}
\DeclareMathOperator{\Aut}{Aut}
\DeclareMathOperator{\Rm}{Rm}
\DeclareMathOperator{\Vol}{Vol}
\DeclareMathOperator{\tr}{tr}
\DeclareMathOperator{\diag}{diag}
\newcommand{\id}{\mathrm{id}}
\newcommand{\st}{\operatorname{st}}
\allowdisplaybreaks[1]
\setlist[enumerate]{label=(\arabic*),leftmargin=2.1em}
\hypersetup{pdftitle={Simplex volumes and simplicial volume of convex projective manifolds},
 pdfsubject={Research manuscript},
 pdfkeywords={convex projective manifolds, simplicial volume, affine spheres, Holmes--Thompson volume}}
\title[Simplicial volumes of convex projective manifolds]{Positivity of simplicial volume\\of convex projective manifolds}
\author{Andrea Tamburelli}
\begin{document}
\begin{abstract}
Let $\Omega$ be an irreducible properly convex projective domain admitting
a discrete cocompact group action of projective automorphisms. We prove that the
Holmes--Thompson volumes of projective simplices with vertices in
$\Omega$ are uniformly bounded. As a consequence, we prove that a closed properly convex projective manifold of dimension at least two has positive simplicial volume if and only if it is irreducible. 
\end{abstract}
\maketitle
\setcounter{tocdepth}{1}
\tableofcontents

\section*{Introduction}\label{sec:introduction}

One of the guiding problems in the study of properly convex projective
manifolds is to understand how the geometry of their developing domains
constrains the topology of compact quotients.  In the divisible strictly convex
case, the Hilbert geometry is Gromov hyperbolic and many features of
negatively curved manifolds persist (\cite{Benoist08Survey}).  General divisible properly convex
domains, however, can have nontrivial boundary segments and properly
embedded simplices.  In particular, neither negative sectional
curvature nor hyperbolicity of the fundamental group is available in
general.  

Recall that a properly convex domain $\Omega\subset\RP^n$ is an open
convex subset that is bounded in some affine chart.  It is
\emph{divisible} if a discrete subgroup $\Gamma<\Aut(\Omega)$ acts
cocompactly, and a closed properly convex projective manifold is a
quotient $M=\Gamma\backslash\Omega$ with $\Gamma$ torsion-free.  The
Hilbert metric on $\Omega$ is a complete projectively invariant
Finsler metric.  Its geodesics include the ordinary projective
segments, although uniqueness of geodesics can fail when the boundary
contains segments.

For an oriented closed $n$-manifold, its simplicial volume is the
$\ell^1$-seminorm of the real fundamental class,
\[
 \|M\|=\inf\left\{\sum_j|a_j| \ \Big| \
       \sum_j a_j\sigma_j\text{ represents }[M]\right\}.
\]
In hyperbolic geometry, positivity of $\|M\|$ follows by integrating
the volume form over straight simplices, whose volumes have a
uniform upper bound (\cite{HM}).  This produces a
bounded cocycle with positive evaluation on the fundamental class,
which is the geometric paradigm underlying Gromov's theory
(\cite{Gromov82}).  We establish an analogous bound for projective simplices in $\Omega$, measured using the
Holmes--Thompson volume of the Hilbert metric.

Let $C\subset\R^{n+1}$ be either of the two proper open cones projecting
to $\Omega$.  We say that $\Omega$ is \emph{irreducible} when there are
no nonzero complementary subspaces $L_1,L_2$ and proper open convex
cones $C_i\subset L_i$ with $C=C_1\times C_2$.  Equivalently,
$\Omega$ is not a nontrivial projective join.  This is the same
irreducibility convention used in the rigidity theory of divisible
properly convex domains (also called \emph{indecomposability} in the
sense of Vey and Benoist); it should not be confused with any
connected-sum decomposition of $M$.  For a divisible domain, Vey's
irreducibility theorem says that this condition implies that every
finite-index subgroup of a dividing group $\Gamma$ acts irreducibly on
$\R^{n+1}$; see \cite[Theorem~5.1]{Benoist08Survey}.  

\begin{mainresult}\label{thm:main}
Let $\Omega\subset\RP^n$, $n\geq2$, be an irreducible divisible
properly convex domain.  There exists $V_\Omega<\infty$ such that
\[
 \Vol_{\mathrm{HT},\Omega}(\sigma)\leq V_\Omega
\]
for every (possibly degenerate) projective $n$-simplex $\sigma$ with
vertices in $\Omega$.  The volume is the Holmes--Thompson volume of
the Hilbert metric of $\Omega$.
\end{mainresult}

The assertion concerns simplices in the universal cover.  In
particular, their vertices can approach arbitrary faces of
$\partial\Omega$ and their projections to a compact quotient need not
be embedded.  The natural
Blaschke volume of the canonical affine sphere is comparable with Holmes--Thompson volume on a compact quotient.  Equivariant
straightening of singular simplices therefore gives the principal
topological applications.

\begin{maincorollary}\label{cor:positive}
Every closed properly convex real projective manifold of dimension
$n\geq2$ whose projective domain is irreducible has positive
simplicial volume.
\end{maincorollary}

\begin{maincorollary}\label{cor:alternative}
Let $M=\Gamma\backslash\Omega$ be a closed properly convex real
projective manifold of dimension $n\geq2$.  Then $\|M\|>0$
if and only if $\Omega$ is irreducible.
\end{maincorollary}

The main ingredient towards the proof of Theorem \ref{thm:main} is a characterization of irreducible divisible domains $\Omega$ in terms of the geometry of the associated affine sphere asymptotic to $C$. This carries the Blaschke metric $h$, which Calabi proved to have nonpositive Ricci tensor \cite[Theorem~5.1]{Calabi72}. We give a proof of this inequality in the cocompact setting needed here in
Proposition~\ref{prop:calabi} and improve this estimate in the case of existence of directions annihilating the Ricci tensor.

\begin{mainresult}\label{thm:nullcase}
Let $\Omega\subset\RP^n$, $n\geq2$, be a divisible properly convex
domain, $C$ its positive cone, and $h$ the Blaschke metric of the
canonical affine sphere, normalized to affine mean curvature $-1$.
If
\[
 \operatorname{Ric}_{h,p}(v,v)=0
 \qquad\text{for some }p\in\Omega\text{ and }0\ne v\in T_p\Omega,
\]
then $C$ is a nontrivial direct product of proper open convex cones.
Consequently, if $C$ is irreducible, there is $\varepsilon>0$ such
that $\operatorname{Ric}_h\leq-\varepsilon h$ throughout $\Omega$.
\end{mainresult}

The two-dimensional equality case is known: see for example
Li--Li--Simon \cite{LLS04}. A similar argument appears in 
Hu--Li--Vrancken \cite{HuLiVrancken08} who give characterizations of Calabi products under algebraic conditions on the cubic tensor. Section~\ref{sec:ricci} shows that the null
distribution is parallel and produces an ambient one-parameter
automorphism whose eigenspaces split $C$.

\subsection*{Related work}

The relation between the Euler characteristic and simplicial volume
has recently advanced substantially.  In \cite{BucherFrancavigliaRago26},
Bucher--Francaviglia--Rago 
construct bounded angle cocycles using immersed straightenings.  In
particular, their Theorem~1.2 proves that for any closed properly
convex real projective $d$-manifold there is a constant $c_d>0$,
with $c_d\leq(d+2)/4$, such that
\begin{equation}\label{eq:BFR}
 |\chi(M)|\leq c_d\|M\|.
\end{equation}
Thus Gromov's question whether vanishing simplicial volume forces
vanishing Euler characteristic has an affirmative answer within the
class of closed properly convex projective manifolds.  This result
should not be conflated with positivity of simplicial volume itself:
\eqref{eq:BFR} gives no information when $\chi(M)=0$.  Indeed,
Bucher--Francaviglia--Rago explicitly ask in their Question~1.4
whether closed convex real projective manifolds with
irreducible universal cover necessarily have positive
simplicial volume.  Corollary~\ref{cor:positive} is an affirmative
answer to precisely that question. Our method is different from the angle-cocycle approach.  Rather than representing the Euler class by a bounded cocycle, we produce a uniform volume bound for projective simplices and hence a bounded
representative of the \emph{volume class}.  \\

Several previously known positivity results provide useful points
of comparison.  For an aspherical closed manifold with word-hyperbolic
fundamental group, positivity follows from bounded-cohomological
results of Mineyev \cite{Mineyev01}.  For closed locally symmetric
spaces of noncompact type, positivity was established by
Lafont--Schmidt \cite{LafontSchmidt06}, with the relevant exceptional
case treated by Bucher--Karlsson \cite{Bucher07}.  These cover
strictly convex divisible domains and irreducible symmetric divisible domains,
respectively, but not the full range of irreducible properly convex
domains containing boundary segments. 

Positivity is also known for closed aspherical manifolds of dimension
at least two whose fundamental groups are non-elementary relatively
hyperbolic. This follows from the relative bounded cohomology theorem
of Mineyev--Yaman, for which Franceschini gives another proof
\cite[Theorem~1.1]{Fra18}. In convex projective geometry, this applies
to certain non-strictly convex divisible domains whose dividing groups
are relatively hyperbolic with respect to stabilizers of maximal
properly embedded simplices, as well as to examples with peripheral
subgroups arising from properly embedded cones; see
\cite[Section~1]{BucherFrancavigliaRago26} and the references therein.
Corollary \ref{cor:positive} establishes positivity for all irreducible divisible domains, without a relative hyperbolicity hypothesis.


There are also geometric examples in dimension four worth
emphasizing.  Riolo--Seppi--Slavich
\cite{RioloSeppiSlavich26} construct closed indecomposable convex
projective four-manifolds with positive Euler characteristic by
gluing hyperbolic pieces.  Their simplicial volume is known to be
positive by the corresponding decomposition results, as noted in
\cite{BucherFrancavigliaRago26}.  Our theorem is
not restricted to such decompositions.

\subsection*{Strategy of the proof}
An elementary example illustrates the obstruction that our
proof must exclude.  Consider the open projective triangle
\[
 \Omega_{\mathrm{tri}}
 =\{[x:y:1]\in\RP^2:x,y>0\}.
\]
Positive dilations of the two affine coordinates are projective
automorphisms of this domain.  Invariance under these dilations
implies that its Holmes--Thompson measure has the form
\[
 d\Vol_{\mathrm{HT},\Omega_{\mathrm{tri}}}
 =c_{\mathrm{HT}}\frac{dx\,dy}{xy},
 \qquad c_{\mathrm{HT}}>0.
\]
In particular, this measure is a constant multiple of Euclidean
area in the coordinates $(\log x,\log y)$.
The transformations $(x,y)\mapsto(ex,y)$ and
$(x,y)\mapsto(x,ey)$ generate a group dividing
$\Omega_{\mathrm{tri}}$: in logarithmic coordinates, their
action is simply translation by $\mathbb Z^2$.
Let $\gamma(x,y)=(ex,y)$, and consider the projective triangles
\[
 \sigma_m=\operatorname{conv}
 \bigl\{(1,1),(1,3),(e^m,1)\bigr\},
 \qquad m\geq1.
\]
Here we keep the first two vertices fixed and move the third
along the orbit $\gamma^m(1,1)$.  Each triangle contains the
rectangle
\[
 \left[1,\frac{1+e^m}{2}\right]\times[1,2].
\]
Consequently,
\[
 \Vol_{\mathrm{HT},\Omega_{\mathrm{tri}}}(\sigma_m)
 \geq c_{\mathrm{HT}}\log2\,
       \log\left(\frac{1+e^m}{2}\right)
 \longrightarrow\infty.
\]
Geometrically, these triangles contain a strip whose width
in logarithmic coordinates stays fixed while its length
tends to infinity. 

\noindent The singular values of the transformations producing this
growth have a particularly simple form.  The determinant-one representative of $\gamma^m$ is
\[
 \operatorname{diag}
 \left(e^{2m/3},e^{-m/3},e^{-m/3}\right).
\]
There are only two distinct singular values: one direction
expands relative to the other two, which are rescaled equally. 
Thus the volume growth occurs along a sequence for which
the total distortion, i.e. the difference between the logarithms of the top and bottom singular values, grows without bound, while the variation
within each of the two groups of singular values remains zero.
Our proof obtains a quantitative obstruction to this behaviour
for irreducible divisible domains and uses it to control
the volumes of simplices. \\

We now give some details on how this is achieved. 
First, we strengthen Calabi's nonpositivity of the Blaschke Ricci tensor:
we show that if the Ricci tensor vanishes on a single nonzero vector then the cone splits as a nontrivial product.  Irreducibility excludes this
possibility, and cocompactness then gives a uniform estimate
$\operatorname{Ric}_h\leq-\varepsilon h$ for some $\varepsilon>0$.

The Gauss map of the hyperbolic affine sphere converts this curvature bound into
a restriction on the singular values of cone automorphisms.
Put $N=n+1$, and let
$\mu_1(g)\geq\cdots\geq\mu_N(g)$ be the logarithms of the singular
values of a representative normalized by $|\det g|=1$.  Define
\[
 \begin{aligned}
 R(\mu(g))&=\mu_1(g)-\mu_N(g),\\
 D_k(\mu(g))&=\sum_{i\leq k}(\mu_1(g)-\mu_i(g))
          +\sum_{i>k}(\mu_i(g)-\mu_N(g)),
          \qquad 1\leq k<N.
 \end{aligned}
\]
Here $R$ measures the total spread of the logarithmic singular
values, while $D_k$ measures their variation within the first
$k$ entries and within the remaining entries.  In particular,
$D_k$ vanishes exactly when the entries are constant within
each of these two groups.  Such patterns arise naturally on
a product cone, where one can rescale the two factors
independently. In the example above we had $R(\mu(\gamma^m))=m$ and $D_1(\mu(\gamma^m))=0$, thus those are the quantities that we need to control in order to exclude the phenomenon described in the previous example. We achieve this by proving the quantitative estimate 
\begin{equation}\label{eq:introGap}
 D_k(\mu(\gamma))
 \geq \delta R(\mu(\gamma))-c_{\mathrm{gap}},
 \qquad
 \forall\gamma\in\Gamma,\quad 1\leq k<N,
\end{equation}
for some $\delta>0$. This is a consequence of the relation between the Ricci curvature and the displacement of a fixed point in the symmetric space of positive definite quadratic forms. The estimate in particular says that as the distortion of a group element grows, its singular values cannot become concentrated in two groups with nearly equal values within each group.

As illustrated in the example above, this controls the volume of simplices. Consider
a projective image $\mathcal D$ of $\Omega$ contained in the
positive orthant, equipped with the measure
\[
 \frac{dy_1\cdots dy_n}{y_1\cdots y_n}
 =d(\log y_1)\cdots d(\log y_n).
\]
Our intermediate goal is to bound its total mass uniformly
over all such projective images.  In logarithmic coordinates,
this is an ordinary Euclidean volume problem. Suppose that $\mathcal D$ contains boxes centred at $a$ and $b$,
with half-widths $r a_i$ and $r b_i$ in the $i$-th coordinate,
where $0<r\leq1/4$.  Taking convex combinations of these boxes
shows that, for every $0\leq\lambda\leq1$, the domain $\mathcal D$ contains
the box centred at
\[
 p^\lambda=(1-\lambda)a+\lambda b
\]
with half-widths $r p_i^\lambda$.  Dividing the $i$-th coordinate
by $p_i^\lambda$ sends this box to the fixed box $[1-r,1+r]^n$.
We compare these coordinate changes as the centre moves along
the segment.
Write $\lambda=e^{-s}$ and put $\xi_i=\log(b_i/a_i)$.  Then
\[
 p_i(s)=(1-e^{-s})a_i+e^{-s}b_i.
\]
For $s\geq\log2$, the factor $1-e^{-s}$ lies between $1/2$ and $1$.
Thus the size of $p_i(s)$ is determined, up to a fixed factor,
by the larger of $a_i$ and $e^{-s}b_i$.  Equivalently,
\[
 \log p_i(s)=
 \begin{cases}
  \log b_i-s+O(1),& s<\xi_i,\\
  \log a_i+O(1),& s\geq\xi_i,
 \end{cases}
\]
where the errors are bounded by $\log2$.  Each coordinate
therefore has a simple model: its logarithm decreases with
slope $-1$ until $s=\xi_i$, and is constant afterwards.

Consider an interval $[s_1,s_2]$ containing no transition value
$\xi_i$ in its interior, and set $\tau=s_2-s_1$.  Suppose that
$k$ coordinates are still in the decreasing part of this model.
For those coordinates,
$p_i(s_1)/p_i(s_2)=e^{\tau+O(1)}$; for the others, this ratio
is $e^{O(1)}$.  The transformation comparing the two
normalizations is represented in homogeneous coordinates by
\[
 \operatorname{diag}\left(
 \frac{p_1(s_1)}{p_1(s_2)},\ldots,
 \frac{p_n(s_1)}{p_n(s_2)},1\right).
\]
Its logarithmic singular values therefore consist of $k$
entries close to $\tau$ and $N-k$ entries close to zero.
Subtracting a common constant to normalize the determinant
does not change $R$ or $D_k$.  Consequently, this transformation
has $R=\tau+O(1)$ and $D_k=O(1)$.

We next relate these coordinate changes to elements of
$\Gamma$.  After composing all the normalizations with one
fixed projective map taking the positive orthant to a bounded
simplex, the normalized domains lie in a fixed bounded set
and contain a ball of radius comparable to $r$ around a fixed
point.  Cocompactness allows us to move the preimage of that
point into a fixed compact subset of $\Omega$.  An elementary
matrix estimate implies that on an interval $[s_1,s_2]$ as above this group element 
$\gamma\in\Gamma$ satisfies
\[
 \begin{aligned}
 R(\mu(\gamma))&=\tau+O(1+\log(1/r)),\\
 D_k(\mu(\gamma))&=O(1+\log(1/r)).
 \end{aligned}
\]
All errors are independent of $\tau$ and of the box centres.
Inequality~\eqref{eq:introGap} therefore bounds $\tau$ by
a constant times $1+\log(1/r)$.

Starting at $s=\log2$ and ending at the largest transition
value, there are at most $n$ such intervals with at least
one changing coordinate.  Summing their lengths bounds
$\max_i\log(b_i/a_i)$; if no transition value exceeds $\log2$,
this upper bound is immediate.  Interchanging $a$ and $b$
gives the corresponding bound for the opposite ratios.
We conclude that
\[
 \max_i|\log(b_i/a_i)|
 \leq c_{\mathrm{sep}}\bigl(1+\log(1/r)\bigr).
\]
Thus boxes of a fixed relative thickness cannot have centres
arbitrarily far apart in logarithmic coordinates.

The domain may still extend far by becoming increasingly thin.
To control this possibility, we partition the orthant into
dyadic boxes.  A cell carrying mass at least $\rho$ contains
a coordinate box with relative widths comparable to $\rho$.
The separation estimate therefore bounds the number of such
cells by a constant times $(1+|\log\rho|)^n$.  Integrating this
counting bound over $\rho$ gives a uniform bound for the total
mass of $\mathcal D$.

The singular-value estimate also holds for the dual action.
We can therefore apply the same integral bound to projective
images of $\Omega^*$.  Faifman's projective duality
\cite[Corollary~3.11]{Faifman21} converts it into a uniform
Funk volume bound for simplices, which we then promote to a bound on the Holmes--Thompson volume.
Finally, equivariant straightening and integration of the
Blaschke volume form produce a bounded cocycle whose evaluation
on the fundamental class is positive.  This proves positivity
of simplicial volume.

\subsection*{Organization of the paper} 

Section~\ref{sec:prelim} fixes the cone, affine-sphere, singular-value,
and Holmes--Thompson conventions.  Section~\ref{sec:ricci} proves Theorem \ref{thm:nullcase}.  Section~\ref{sec:cartan} establishes
the Cartan gap; Section~\ref{sec:packing} proves the integral bound estimate.  Section~\ref{sec:volumes} applies duality to obtain
the uniform bound in Theorem \ref{thm:main}.  Section~\ref{sec:cocycle} concludes with an explicit equivariant straightening and the proof of the
corollaries regarding the simplicial volume.

\subsection*{AI declaration} The main strategy of the proof was motivated by the example, explained in the introduction, of a family of simplices with unbounded Holmes--Thomson volume, which was an original observation of the author, who thus realized that the Cartan gap \eqref{eq:introGap} could be the key to prove the main theorem. The Ricci estimate and the Cartan gap were proved by the author without any AI tool. Chat GPT Sol 6.1 was later used to deduce the integral bound in Section \ref{sec:packing} from the Cartan gap. In other words, all mathematical ideas outside Section \ref{sec:packing} are original of the author. The author also used Chat GPT Sol 6.1 for literature review and editing.

\subsection*{Acknowledgment} The author acknowledges support from the Italian Ministry of University and Research (MUR) through the FIS3 project “Geometry and Topology of higher Teichm\"uller spaces”, grant no. FIS-2024-05142.

\section{Preliminaries}\label{sec:prelim}

We write $n=\dim\Omega\geq2$ and $N=n+1$. We fix the standard
ambient volume density on $\R^N$.
The standard Euclidean operator norm of a linear map is denoted by
$\|\cdot\|$. 

\subsection{Cones, duality, and canonical lifts}

An open convex cone $C\subset\R^N$ is \emph{proper} if it is nonempty
and full-dimensional and its closure is pointed:
$\overline C\cap(-\overline C)=\{0\}$. Its projectivization is a
properly convex domain. Conversely, a properly convex domain has two
opposite proper cones; we fix one and call it $C$. Define its open dual
cone and projective dual domain by
\[
 C^*=\{\ell\in(\R^N)^*: \ell(v)>0
                  \text{ for every }v\in\overline C\setminus\{0\}\},
 \qquad \Omega^*=\mathbb P(C^*).
\]
In particular, an element of $C^*$ is strictly positive on
every nonzero boundary ray of $C$.

A product decomposition of $C$ consists of a direct sum
$\R^N=L_1\oplus L_2$ and proper open convex cones $C_i\subset L_i$
such that $C=C_1\times C_2$. It is nontrivial when both $L_i$ are
nonzero; a one-dimensional factor is allowed. Projectivizing gives the
join of $\mathbb P(C_1)$ and $\mathbb P(C_2)$. 

Every $g\in\Aut(\Omega)$ has a unique linear representative preserving
$C$ and satisfying $|\det g|=1$ (\cite{BallasDanciger}). We use the same letter for a projective automorphism and this representative. The contragredient group acts
on $C^*$ by $g^{-\mathsf T}$.

\subsection{The canonical affine sphere}

We use the existence and uniqueness theorem for complete hyperbolic
affine spheres, together with its conormal duality; see
\cite[Theorem~3 and Proposition~8]{Loftin08}. We record the exact
normalizations needed here.

\begin{lemma}\label{lem:canonical}
There is a unique radial lift $f\colon\Omega\to C$, with
$[f(x)]=x$, whose image is a hyperbolic affine sphere centered at the
origin, of affine mean curvature $-1$, meeting every ray of $C$ once
and having complete Blaschke metric $h$. It is equivariant for the cone-preserving lifts just
defined. Its conormal map $\nu\colon\Omega\to C^*$, characterized by
\[
 \nu_x(f(x))=1,\qquad \nu_x(f_*X)=0,
\]
projectivizes to a diffeomorphism $\Omega\to\Omega^*$ and is
equivariant for the contragredient action. In particular, if $\Gamma$
divides $\Omega$, its contragredient group divides $\Omega^*$.
\end{lemma}

Let $D$ be the flat ambient connection and let $\nabla$ be the
Levi--Civita connection of $h$. If $\nabla^{\mathrm{aff}}$ is the
induced affine connection, set $A=\nabla^{\mathrm{aff}}-\nabla$. The
structure equation is
\begin{equation}\label{eq:structure}
 D_X f_*Y=f_*(\nabla_XY+A_XY)+h(X,Y)f.
\end{equation}
The tensor $A(X,Y,Z)=h(A_XY,Z)$ is symmetric in all three entries,
each endomorphism $A_X$ is trace free, and
\begin{equation}\label{eq:codazzi}
 (\nabla_XA)(Y,Z,W)
 \quad\text{is symmetric in }X,Y,Z,W.
\end{equation}
We use the conventions
\[
 \begin{split}
 \Rm(X,Y)Z&=\nabla_X\nabla_YZ-\nabla_Y\nabla_XZ-\nabla_{[X,Y]}Z,\\
 \operatorname{Ric}(Y,Z)&=\sum_a h(\Rm(e_a,Y)Z,e_a),\qquad
 \Delta=\tr_h\nabla^2.
 \end{split}
\]
Thus $\Delta$ has positive principal part $\sum_a\partial_a^2$ in a
normal frame. The Gauss equation, with these conventions, is
\begin{equation}\label{eq:gauss}
 \Rm(X,Y)Z=-h(Y,Z)X+h(X,Z)Y-[A_X,A_Y]Z.
\end{equation}
In an $h$-orthonormal frame put $A_i=A_{e_i}$ and
$B_{ij}=\tr(A_iA_j)$. Tracing~\eqref{eq:gauss} gives
\begin{equation}\label{eq:ricci}
 \operatorname{Ric}_h=-(n-1)h+B.
\end{equation}
The tensor $B$ is symmetric and nonnegative. We freely
identify symmetric two-tensors with self-adjoint endomorphisms using
$h$.

\subsection{Cartan projections}

For a linear map $g$ with $|\det g|=1$, let
$s_1(g)\geq\cdots\geq s_N(g)>0$ be its singular values. Its
\emph{Cartan vector} is
\[
 \mu(g)=(\mu_1(g),\ldots,\mu_N(g)),\qquad
 \mu_i(g)=\log s_i(g).
\]
It lies in the closed ordered chamber of the trace-zero hyperplane.
For an arbitrary projective transformation, we use any representative
normalized to absolute determinant one; the two sign choices have the
same singular values. On an ordered trace-zero vector $\mu$ set
\begin{equation}\label{def:cartan}
 \begin{split}
 R(\mu)&=\mu_1-\mu_N,\\
 D_k(\mu)&=\sum_{i\leq k}(\mu_1-\mu_i)
              +\sum_{i>k}(\mu_i-\mu_N),\qquad 1\leq k<N.
 \end{split}
\end{equation}
Thus $D_k$ vanishes exactly on the ray of vectors with their first
$k$ entries equal and their remaining entries equal. 

We record the elementary perturbation estimate used below.

\begin{lemma}\label{lem:cartanperturb}
Suppose $g=\mathsf U g_0\mathsf V$, where $g,\mathsf U,g_0,\mathsf V$ have determinant of absolute value one. If
\[
 \|\mathsf U\|,\|\mathsf U^{-1}\|\leq e^{\eta_{\mathsf U}},\qquad
 \|\mathsf V\|,\|\mathsf V^{-1}\|\leq e^{\eta_{\mathsf V}},
\]
then $|\mu_i(g)-\mu_i(g_0)|\leq \eta_{\mathsf U}+\eta_{\mathsf V}$ for every $i$. In
particular, a fixed change of ambient basis changes all Cartan
coordinates by a uniformly bounded amount. Moreover, the dual Cartan
vector satisfies
\begin{equation}\label{eq:dualcartan}
 \widetilde\mu_i=\mu_i(g^{-\mathsf T})=-\mu_{N+1-i}(g),\qquad
 R(\widetilde\mu)=R(\mu),\qquad
 D_k(\widetilde\mu)=D_{N-k}(\mu).
\end{equation}
\end{lemma}

\begin{proof}
The min--max characterization of singular values gives
$s_i(\mathsf U g_0\mathsf V)\leq\|\mathsf U\|\|\mathsf V\|s_i(g_0)$. Applying the same inequality to
$g_0=\mathsf U^{-1}g\mathsf V^{-1}$ gives the reverse bound. Taking logarithms proves
the first assertion. A basis change is conjugation by a fixed matrix;
normalizing that matrix to absolute determinant one leaves the
conjugation unchanged and gives the second assertion. Inversion
reverses and reciprocates the singular values, whereas transposition
does not change them. This proves the first equality
in~\eqref{eq:dualcartan}; substitution into~\eqref{def:cartan} proves
the other two.
\end{proof}

If two ordered trace-zero vectors differ by at most $e_0$ in every
coordinate, their $R$ values differ by at most $2e_0$, and their $D_k$
values differ by at most $2Ne_0$. These bounds also follow immediately
from~\eqref{def:cartan}.

\subsection{Funk and Holmes--Thompson volumes} 

Choose a bounded affine realization of a properly convex domain $U$.
For $x\in U$ and $v\neq0$, the Funk norm is defined by
$F_U(x,v)=1/t$, where $t>0$ is the unique parameter such that
$x+tv\in\partial U$; put $F_U(x,0)=0$. Its tangent unit ball is
$U-x$. Its cotangent unit ball is the polar convex body
\[
 K^*_{U,x}=\{\xi \ | \ \xi(z-x)\leq1\text{ for every }z\in U\}.
\]
The Hilbert norm is
$H_U(x,v)=\tfrac12(F_U(x,v)+F_U(x,-v))$. In terms of cotangent
unit balls,
\begin{equation}\label{eq:hilbertfunk}
 K^*_{H_U,x}=\tfrac12(K^*_{U,x}-K^*_{U,x}).
\end{equation}
This follows because a norm is the support function of its cotangent
unit ball and support functions add under Minkowski addition.

Let $\omega_n$ be the Euclidean volume of the unit ball in $\R^n$.
For a possibly asymmetric Finsler metric $\mathfrak f$, its Holmes--Thompson
volume density is
\begin{equation}\label{eq:htdensity}
 d\Vol_{\mathrm{HT},\mathfrak f}(x)
       =\frac{\Vol_{\mathrm{Eucl}}(K^*_{\mathfrak f,x})}{\omega_n}\,
                         dx_1\cdots dx_n.
\end{equation}
We abbreviate the Hilbert
Holmes--Thompson measure by $\Vol_{\mathrm{HT},U}$ and the Funk
Holmes--Thompson measure by $\Vol_{F,U}$. Unlike the Hilbert norm,
the Funk norm depends on the affine chart. Its total volume on a
convex subdomain has the projective invariance and duality stated
precisely in Section~\ref{sec:volumes}.

For later use, both Hilbert and Funk cotangent unit balls vary
continuously on compact subsets of a bounded convex domain, in the
Hausdorff topology. One can see this by first translating the closed
domain by $-x$ and then taking polars; on a compact set of basepoints
these bodies contain a fixed ball and lie in a fixed ball, so polarity
is continuous. Their Euclidean volumes are therefore positive and
continuous. In particular the Holmes--Thompson density is a
positive continuous density, invariant under projective
automorphisms.

\section{Ricci curvature and cone splitting}\label{sec:ricci}

Calabi's nonpositivity theorem~\cite[Theorem~5.1]{Calabi72} applies to
complete hyperbolic affine spheres.  We reprove it in
the present cocompact setting, using a differential identity for the
negative of the Ricci tensor. We then analyze the null case: this
identity and a strong minimum principle show that the Ricci kernel has
constant rank and is parallel. A nonzero null direction consequently
produces a parallel field and a nontrivial product decomposition of the
cone, proving Theorem~\ref{thm:nullcase}. \\

We identify symmetric covariant two-tensors with self-adjoint
endomorphisms by means of $h$.  At a point, choose a normal
$h$-orthonormal frame $e_1,\ldots,e_n$, and write
\[
 A_i=A_{e_i},\qquad
 A_{ijk}=h(A_{e_i}e_j,e_k),\qquad
 B_{ij}=\tr(A_iA_j),\qquad P=(n-1)\id-B.
\]
Thus $P=-\operatorname{Ric}_h$ by~\eqref{eq:ricci}.  Repeated indices
are summed from $1$ to $n$.  The inner product of endomorphisms is the Frobenius inner product
$\langle U,V\rangle_{\mathrm{HS}}=\tr(U^*V)$, where the adjoint is taken with
respect to $h$.  Define the symmetric nonnegative two-tensors
\begin{align*}
 \mathcal H_{ij}
  &=\sum_{k,a,b}(\nabla_kA)_{iab}(\nabla_kA)_{jab},\\
 \mathcal Q_{ij}
  &=\sum_a\langle[A_i,A_a],[A_j,A_a]\rangle_{\mathrm{HS}},\\
 \mathcal R^2_{ij}
  &=\sum_a\langle \Rm(e_i,e_a),\Rm(e_j,e_a)\rangle_{\mathrm{HS}}.
\end{align*}

The tensors $A$, $B$ and $P$ satisfy the following differential equations. These computations are classical and can be found for instance in \cite{Hildebrand13} and
\cite[Proposition~3.4]{KlartagKolesnikov17}. We re-derive all three
formulas here in order to avoid painful book-keeping of the different conventions used in the aforementioned references.

\begin{lemma}\label{lem:simons}
The cubic form satisfies
\begin{equation}\label{eq:simons}
\begin{split}
 \Delta A_{ijk}={}&-(n+1)A_{ijk}
   +B_{ip}A_{pjk}+B_{jp}A_{ipk}+B_{kp}A_{ijp}\\
   &-2\tr(A_iA_jA_k).
\end{split}
\end{equation}
Its quadratic contraction satisfies
\begin{equation}\label{eq:deltaB}
 \frac12\Delta B=\mathcal H-(n+1)B+B^2+\mathcal Q,
\end{equation}
and
\begin{equation}\label{eq:deltaP}
 \frac12\Delta P=(n+1)P-P^2-\mathcal R^2-\mathcal H.
\end{equation}
\end{lemma}

\begin{proof}
Since $A$ is tracefree and $\nabla A$ is full symmetric,
we have
\[
 \sum_a\nabla_a A_{ajk}
   =\nabla_j\sum_a A_{aak}=0.
\]
Writing $\Rm_{aij}{}^p$ for the coefficients of
$\Rm(e_a,e_i)e_j$, and using the preceding divergence identity, we obtain
\begin{align*}
 \Delta A_{ijk}
  &=\sum_a(\nabla_a\nabla_i-\nabla_i\nabla_a)A_{ajk}\\
  &=-\sum_{a,p}\bigl(\Rm_{aia}{}^pA_{pjk}
       +\Rm_{aij}{}^pA_{apk}+\Rm_{aik}{}^pA_{ajp}\bigr)\\
  &=\operatorname{Ric}_{ip}A_{pjk}
       -\sum_{a,p}\bigl(\Rm_{aij}{}^pA_{apk}
                         +\Rm_{aik}{}^pA_{ajp}\bigr).
\end{align*}
The Gauss equation~\eqref{eq:gauss} yields
\begin{align*}
 -\sum_{a,p}\Rm_{aij}{}^pA_{apk}
  &=-A_{ijk}+\sum_{a,p}[A_a,A_i]^p{}_jA_{apk}\\
  &=-A_{ijk}+B_{kp}A_{ijp}-\tr(A_iA_jA_k).
\end{align*}
The other curvature contraction interchanges $j$ and $k$.  The trace of
a product of three symmetric matrices is symmetric in those matrices,
by cyclicity of the trace and transposition.  Substituting
$\operatorname{Ric}_{ip}=-(n-1)\delta_{ip}+B_{ip}$ proves
\eqref{eq:simons}.

\noindent For the quadratic contraction, introduce
\[
 B^{(1)}_{ij}=\tr(BA_iA_j),\qquad
 B^{(2)}_{ij}=\sum_a\tr(A_iA_aA_jA_a).
\]
Both tensors are symmetric.  Cubic symmetry gives
\[
 \left(\sum_a A_a^2\right)_{ij}
   =\sum_{a,p}A_{aip}A_{apj}=B_{ij}.
\]
Since every commutator $[A_i,A_a]$ is skew-adjoint, its expansion gives
\begin{align*}
 \mathcal Q_{ij}
  &=-\sum_a\tr\bigl((A_iA_a-A_aA_i)(A_jA_a-A_aA_j)\bigr)=2B^{(1)}_{ij}-2B^{(2)}_{ij}.
\end{align*}
The cubic term of~\eqref{eq:simons} has contraction
\begin{align*}
 \sum_{a,b}A_{jab}\tr(A_iA_aA_b)
  &=\sum_{a,b,p,q,r}A_{jab}A_{ipq}A_{aqr}A_{brp}\\
  &=\sum_r\tr(A_jA_rA_iA_r)=B^{(2)}_{ij}.
\end{align*}
Here the second equality uses $A_{aqr}=A_{raq}$ and
$A_{brp}=A_{rbp}$.  The product rule therefore gives
\begin{align*}
 \frac12\Delta B_{ij}
  &=\mathcal H_{ij}+\frac12\sum_{a,b}
        \bigl((\Delta A_{iab})A_{jab}+A_{iab}(\Delta A_{jab})\bigr)\\
  &=\mathcal H_{ij}-(n+1)B_{ij}+(B^2)_{ij}+2B^{(1)}_{ij}-2B^{(2)}_{ij}.
\end{align*}
This proves~\eqref{eq:deltaB}.

\noindent We next contract the Gauss equation.  With
$J_{ij}X=h(e_j,X)e_i-h(e_i,X)e_j$, one has
\[
 \sum_a\langle J_{ia},J_{ja}\rangle_{\mathrm{HS}}
   =2(n-1)\delta_{ij}.
\]
The mixed contraction is
\begin{align*}
 \sum_a\langle J_{ia},[A_j,A_a]\rangle_{\mathrm{HS}}
  &=2\sum_a[A_j,A_a]_{ia}\\
  &=2\sum_{a,p}\bigl(A_{jip}A_{apa}-A_{aip}A_{jpa}\bigr)
   =-2B_{ij}.
\end{align*}
The first term vanishes by trace freeness.  Consequently
\begin{equation}\label{eq:curvature-square}
 \mathcal R^2=\mathcal Q+2(n-1)\id-4B.
\end{equation}
Combining this with~\eqref{eq:deltaB}, and substituting
$B=(n-1)\id-P$, proves~\eqref{eq:deltaP}.
\end{proof}

\begin{proposition}[Calabi]\label{prop:calabi}
The Blaschke Ricci tensor satisfies
\[
 -(n-1)h\leq\operatorname{Ric}_h\leq 0.
\]
Equivalently, $0\leq B\leq(n-1)\id$ and $P\geq0$.
\end{proposition}

\begin{proof}
Equations~\eqref{eq:deltaB} and~\eqref{eq:curvature-square} give
\begin{equation}\label{eq:calabi-factorization}
 \frac12\Delta B
   =\mathcal H+\mathcal R^2+(B-(n-1)\id)(B+2\id).
\end{equation}
Choose a compact fundamental domain $\mathcal F\subset\Omega$ for the action of $\Gamma$.
The largest eigenvalue of $B$ is continuous and $\Gamma$-invariant,
so its global maximum $\lambda_{\max}$ is attained at a point $p\in \mathcal F$.
Moreover, $\lambda_{\max}\geq0$ because $B$ is symmetric and non-negative definite.

Choose a unit eigenvector $v$ at $p$ relative to $\lambda_{\max}$ and extend it to a
local unit field with $\nabla v(p)=0$.  The scalar function
$b=B(v,v)$ has a local maximum at $p$.  At that point,
\[
 \Delta b=(\Delta B)(v,v)+2B(\Delta v,v)=(\Delta B)(v,v).
\]
The last equality follows from $Bv=\lambda_{\max} v$ and
$h(\Delta v,v)=0$, the latter being the second derivative of the
unit-length condition at a point where $\nabla v=0$.  Hence
\eqref{eq:calabi-factorization} implies
\[
 0\geq\frac12(\Delta B)(v,v)
   =\mathcal H(v,v)+\mathcal R^2(v,v)
      +(\lambda_{\max}-(n-1))(\lambda_{\max}+2).
\]
Both $\mathcal H(v,v)$ and $\mathcal R^2(v,v)$ are nonnegative and $\lambda_{\max}+2>0$, so
$\lambda_{\max}\leq n-1$.  Thus $P\geq0$.  The remaining inequality
$\operatorname{Ric}_h\geq-(n-1)h$ follows directly from $B\geq0$
and~\eqref{eq:ricci}.
\end{proof}

We shall use the following scalar minimum principle for continuous
functions.  

\begin{lemma}\label{lem:viscmin}
Let $(\mathcal X,h)$ be a connected smooth Riemannian manifold and let $c\geq0$
be a constant.  Suppose that $u\colon \mathcal X\to[0,\infty)$ is continuous
and satisfies
\[
 -\Delta u+cu\geq0
\]
in the following sense: whenever a $C^2$ function $\phi$ touches $u$
from below at $p$, one has $\Delta\phi(p)\leq cu(p)$.
If $u$ vanishes at one point, then $u$ vanishes identically.
\end{lemma}

\begin{proof}
This is a special case of the viscosity strong maximum principle
in~\cite[Theorem 4.7, Example 4.12]{GoffiPediconi}, applied to the Laplacian operator.  We give a direct proof in this case for the convenience of the reader.

Suppose that both the positive set and the zero set of $u$ are
nonempty.  In a coordinate chart choose a Euclidean ball $B_{r_0}(x)$
contained in the positive set and tangent to its zero set at a point
$q$ of the chart.  Such a ball is obtained by taking its center in the
positive set sufficiently close to its boundary and using the distance
to the zero set as radius.  It can be chosen with its closure contained
in the chart.  In these coordinates write
\[
 \Delta=a^{ij}(y)\partial_i\partial_j+b^i(y)\partial_i,
\]
where $(a^{ij})$ is uniformly positive definite on a neighborhood of
$\overline{B_{r_0}(x)}$.  On the annulus $r_0/2\leq|y-x|\leq r_0$, set
\[
 w(y)=e^{-\alpha|y-x|^2}-e^{-\alpha r_0^2}.
\]
For sufficiently large $\alpha$,
\begin{equation}\label{eq:minimum-barrier}
 \Delta w-cw>0
\end{equation}
on a neighborhood of this closed annulus.  Indeed, writing
$\xi=y-x$, the left side is
\[
 e^{-\alpha|\xi|^2}
 \left(4\alpha^2a^{ij}\xi_i\xi_j
       -2\alpha\sum_i a^{ii}-2\alpha b^i\xi_i-c\right)
       +ce^{-\alpha r_0^2};
\]
uniform ellipticity and $|\xi|\geq r_0/2$ make the quadratic term in
$\alpha$ dominate the other terms.

Choose $\delta>0$ so small that $u\geq\delta w$ on the inner sphere.
The same inequality holds on the outer sphere, where $w=0$.  If
$u-\delta w$ had a negative minimum $m$ inside the annulus, then
$\delta w+m$ would touch $u$ from below there.  The assumed viscosity
inequality and~\eqref{eq:minimum-barrier} would give
\[
 0\leq-\delta\Delta w+c(\delta w+m)<0,
\]
a contradiction.  Thus $u\geq\delta w$ on the annulus.  Immediately
outside the outer sphere one has $w<0\leq u$, so $\delta w$ touches
$u$ from below at $q$.  Since $u(q)=w(q)=0$,
\eqref{eq:minimum-barrier} contradicts the viscosity inequality at
$q$ as well.  This proves the lemma.
\end{proof}

\begin{lemma}\label{lem:kernelparallel}
The kernel $\mathcal N=\ker P$ is a
smooth parallel subbundle.  Moreover, for every local section $v$ of $\mathcal N$,
\begin{equation}\label{eq:null-consequences}
 \Rm(v,X)=0,\qquad (\nabla_XA)(v,Y,Z)=0
\end{equation}
for all tangent fields $X,Y,Z$.
\end{lemma}

\begin{proof}
Let $s_k$ be the sum of the $k$ smallest eigenvalues of $P$, for
$1\leq k\leq n$.  These functions are continuous and nonnegative by
Proposition~\ref{prop:calabi}.  We claim that they are viscosity
supersolutions of
\begin{equation}\label{eq:eigenvalue-supersolution}
 -\Delta s_k+2(n+1)s_k\geq0.
\end{equation}
Suppose that $\phi$ touches $s_k$ from below at $p$.  Choose
orthonormal eigenvectors $e_1,\ldots,e_k$ for the $k$ smallest
eigenvalues at $p$, and extend them orthonormally with
$\nabla e_i(p)=0$.  The variational characterization of the sum of
the smallest eigenvalues gives
\[
 s_{\mathrm{fr}}=\sum_{i=1}^kP(e_i,e_i)\geq s_k,\qquad s_{\mathrm{fr}}(p)=s_k(p).
\]
Thus $\phi$ also touches the smooth function $s_{\mathrm{fr}}$ from below at $p$.
All frame-derivative terms in its Laplacian vanish there: the terms
with first derivatives vanish by construction, while
\[
 P(\Delta e_i,e_i)
  =\lambda_i(P)h(\Delta e_i,e_i)=0.
\]
It follows from~\eqref{eq:deltaP} that
\[
 \Delta\phi(p)\leq\Delta s_{\mathrm{fr}}(p)
   =\sum_{i=1}^k(\Delta P)(e_i,e_i)
   \leq2(n+1)s_k(p),
\]
because $P^2$, $\mathcal R^2$, and $\mathcal H$ are nonnegative.  This proves
\eqref{eq:eigenvalue-supersolution} with exactly the lower-test
interpretation in Lemma~\ref{lem:viscmin}.

Let $d_{\mathrm{null}}$ be the largest nullity attained by $P$.  If $d_{\mathrm{null}}=0$, its nullity
is already constant.  Otherwise $s_{d_{\mathrm{null}}}$ vanishes somewhere, and
Lemma~\ref{lem:viscmin} gives $s_{d_{\mathrm{null}}}\equiv0$.  The nullity is therefore
at least $d_{\mathrm{null}}$ everywhere, and maximality makes it exactly $d_{\mathrm{null}}$.  The
constant-rank theorem now makes $\mathcal N=\ker P$ a smooth subbundle.

For a local section $v$ of $\mathcal N$, differentiating $P(v,\cdot)=0$ gives
\[
 (\nabla_aP)(v,W)=-P(\nabla_av,W).
\]
Twice differentiating the identically zero scalar $P(v,v)$, in a
normal orthonormal frame at the point of calculation, gives
\begin{align*}
 0={}&(\Delta P)(v,v)
       +4\sum_a(\nabla_aP)(\nabla_av,v)
       +2P(\Delta v,v)
       +2\sum_aP(\nabla_av,\nabla_av)\\
   ={}&(\Delta P)(v,v)-2\sum_aP(\nabla_av,\nabla_av).
\end{align*}
On the other hand, $Pv=0$ in~\eqref{eq:deltaP} gives the exact identity
\begin{equation}\label{eq:null-sum}
 \sum_aP(\nabla_av,\nabla_av)+\mathcal R^2(v,v)+\mathcal H(v,v)=0.
\end{equation}
All summands are nonnegative.  Hence $\nabla_av\in \mathcal N$, which proves
that $\mathcal N$ is parallel.  The equality $\mathcal R^2(v,v)=0$ means
\[
 \sum_a\lVert \Rm(v,e_a)\rVert_{\mathrm{HS}}^2=0,
\]
and the equality $\mathcal H(v,v)=0$ means
\[
 \sum_{k,a,b}\bigl((\nabla_{e_k}A)(v,e_a,e_b)\bigr)^2=0.
\]
These are precisely~\eqref{eq:null-consequences}.
\end{proof}

\begin{proposition}\label{prop:riccisplit}
If the Blaschke Ricci tensor has a nonzero null vector at one point,
then the cone $C$ is a nontrivial direct product of proper open convex
cones.
\end{proposition}

\begin{proof}
By Lemma~\ref{lem:kernelparallel}, the parallel bundle $\mathcal N=\ker P$ has
positive rank and satisfies $\Rm(v,X)=0$ for $v\in \mathcal N$.  The first Bianchi
identity gives
\[
 \Rm(X,Y)v=-\Rm(Y,v)X-\Rm(v,X)Y=0.
\]
Thus the induced connection on $\mathcal N$ is flat.  Since $\Omega$ is simply
connected, every vector in a fibre of $\mathcal N$ extends to a unique global
parallel section of $\mathcal N$.

Let $\mathscr V$ be the full vector space of global parallel vector fields on
$(\Omega,h)$.  Every such field $W$ belongs to $\mathcal N$: parallelism gives
$\Rm(X,Y)W=0$, and hence $\operatorname{Ric}(W,W)=0$; since $P\geq0$,
this implies $PW=0$.  Consequently evaluation at any point identifies
$\mathscr V$ with the corresponding fibre of $\mathcal N$.  In particular $\mathscr V$ is a
nonzero finite-dimensional Euclidean space, with the inner product
given by $h$ on parallel fields.

For $X,Y\in \mathscr V$, equation~\eqref{eq:null-consequences} gives, for
arbitrary $Z,W$,
\[
 h(\nabla_Z(A_XY),W)=(\nabla_ZA)(X,Y,W)=0.
\]
Therefore $A_XY$ is a global parallel field, proving $A(\mathscr V,\mathscr V)\subset \mathscr V$.
The same calculation shows that $A(X,Y,Z)$ is constant for
$X,Y,Z\in \mathscr V$.  We may thus regard the restriction of $A$ as a fixed
symmetric cubic on $\mathscr V$.

Choose a unit vector $u\in \mathscr V$ at which $A(u,u,u)$ is maximal on the
unit sphere of $\mathscr V$.  The Lagrange multiplier equation, together with
$A_uu\in \mathscr V$, gives
\begin{equation}\label{eq:cubic-maximizer}
 A_uu=a u,\qquad a=A(u,u,u),
\end{equation}
where $a$ is constant.  By self-adjointness $A_u$ preserves $u^\perp$.
For $X\perp u$, the Gauss equation and
\eqref{eq:null-consequences} give
\[
 0=\Rm(u,X)u=X-A_u^2X+aA_uX.
\]
Thus
\begin{equation}\label{eq:Au-polynomial}
 A_u^2-aA_u-\id=0\quad\hbox{on }u^\perp.
\end{equation}
Moreover $u$ and $A_u$ are parallel, the latter because
\[
 h((\nabla_XA_u)Y,Z)=(\nabla_XA)(u,Y,Z)=0.
\]
The radial vector $f$ is transverse to the tangent space of the
affine sphere.  Define an ambient endomorphism $\Theta$ at each point by
\begin{equation}\label{eq:ambient-splitting-map}
 \Theta f=f_*u,\qquad \Theta f_*Y=f_*A_uY+h(u,Y)f.
\end{equation}
We show that it is a constant matrix.  Using the structure
equation~\eqref{eq:structure} and the parallelism just proved gives
\[
 (D_X\Theta)f=f_*(A_Xu-A_uX)=0.
\]
On tangent vectors, the same differentiation gives
\begin{align*}
 (D_X\Theta)f_*Y
  ={}&f_*\bigl([A_X,A_u]Y+h(u,Y)X-h(X,Y)u\bigr)\\
    &+\bigl(h(X,A_uY)-h(u,A_XY)\bigr)f\\
  ={}&-f_*\Rm(X,u)Y=0.
\end{align*}
The radial coefficient vanishes by cubic symmetry, and the last
equality follows from~\eqref{eq:null-consequences}.  These equalities
hold on a basis of the ambient space, so $D_X\Theta=0$ for every $X$.
Connectedness makes $\Theta$ a single constant endomorphism of $\R^N$.
On $\operatorname{span}(f,f_*u)$, its matrix is
\[
 \begin{pmatrix}0&1\\1&a\end{pmatrix}.
\]
The complementary invariant space $f_*(u^\perp)$ has $\Theta=A_u$ and
satisfies~\eqref{eq:Au-polynomial}.  Consequently $\Theta$ is
diagonalizable and has exactly the two distinct real eigenvalues
\[
 \lambda_\pm=\frac{a\pm\sqrt{a^2+4}}2.
\]
Both eigenvalues occur already on the displayed two-dimensional
space.  Denote their nonzero ambient eigenspaces by $L_+$ and $L_-$.

The parallel unit field $u$ is complete, since its integral curves are
unit geodesics of the complete metric $h$.  Let $\phi_t$ be its global
flow.  Equation~\eqref{eq:ambient-splitting-map} gives
\[
 \frac{d}{dt}f(\phi_t(x))=\Theta f(\phi_t(x)),\qquad
 f\circ\phi_t=e^{t\Theta}f.
\]
The affine sphere meets every ray of $C$ once, so $e^{t\Theta}C=C$ for all
$t\in\R$.  Combining this action with radial multiplication allows
independent positive scaling on $L_+$ and $L_-$.  Indeed, prescribed
factors $s_+,s_->0$ are obtained from $r e^{t\Theta}$ by taking
\[
 t=\frac{\log s_+-\log s_-}{\lambda_+-\lambda_-},\qquad
 \log r=\log s_+-t\lambda_+.
\]
Set $C_\pm=\operatorname{pr}_{L_\pm}C$.  These are nonempty open
full-dimensional convex cones in $L_\pm$.  We claim that
\begin{equation}\label{eq:cone-product}
 C=C_+\times C_-
 \quad\hbox{in }\R^N=L_+\oplus L_-.
\end{equation}
For $x_+\in C_+$ and $x_-\in C_-$ choose
$(x_+,y_-)\in C$ and $(y_+,x_-)\in C$.  Independent positive scaling
and addition give
\[
 (x_++\eta y_+,\,x_-+\eta y_-)\in C
 \qquad(\eta>0).
\]
Letting $\eta$ tend to zero proves
$C_+\times C_-\subset\overline C$.  Since the product is open, it is
contained in $\operatorname{int}\overline C=C$.  The reverse inclusion
is immediate from the definitions, proving~\eqref{eq:cone-product}.
Finally
$\overline C=\overline C_+\times\overline C_-$ is pointed, so each
closed factor is pointed.  Thus both open factors are proper.  Their
ambient spaces are nonzero, making this product decomposition
nontrivial.
\end{proof}

\begin{proposition}\label{prop:riccineg}
If $C$ is irreducible, there is a constant $\varepsilon>0$ such that
\begin{equation}\label{eq:ricci-gap}
 -(n-1)h\leq\operatorname{Ric}_h\leq-\varepsilon h
 \quad\hbox{on }\Omega.
\end{equation}
\end{proposition}

\begin{proof}
Proposition~\ref{prop:calabi} gives nonpositive Ricci curvature.
Proposition~\ref{prop:riccisplit} and irreducibility exclude a null
direction at every point, so $P=-\operatorname{Ric}_h$ is positive
definite.  Its smallest eigenvalue is a continuous positive
$\Gamma$-invariant function.  It has a positive minimum
$\varepsilon$ on a compact set $\mathcal F$ with $\Gamma \mathcal F=\Omega$.
This proves the upper bound in~\eqref{eq:ricci-gap}; the lower bound
is Proposition~\ref{prop:calabi}. Together with
Proposition~\ref{prop:riccisplit}, this proves Theorem~\ref{thm:nullcase}.
\end{proof}

\section{The affine-sphere Gauss map and Cartan projections}\label{sec:cartan}

In this section we assume that $C$ is irreducible.  By
Proposition~\ref{prop:riccineg}, there is $\varepsilon>0$ such that
\begin{equation}\label{eq:cartan-ricci-assumption}
 \operatorname{Ric}_h\leq-\varepsilon h
 \qquad\hbox{on }\Omega.
\end{equation}
We write $\nabla$ for the Levi--Civita connection of $h$.  Our aim is to
convert this curvature bound into a restriction on the Cartan
projections of cone automorphisms.  The comparison will use the
Blaschke length of projective segments. 

Fix $p\in\Omega$.  For each $x\in\Omega$, let $G_x$ be the positive
definite form on $\R^N$ characterized by
\[
 G_x(f(x),f(x))=1,\qquad
 G_x(f(x),f_*T_x\Omega)=0,\qquad
 G_x(f_*X,f_*Y)=h_x(X,Y).
\]
We denote the norm defined by $G_p$ by $\|\cdot\|_p$, and its dual norm
by $\|\cdot\|_{p,*}$.  For a cone-preserving matrix $g$ normalized by
$|\det g|=1$, write
\[
 \mu^p_1(g)\geq\cdots\geq\mu^p_N(g),\qquad
 \sum_i\mu^p_i(g)=0,
\]
for its logarithmic singular values in this norm.  We first work with
$\mu^p$ and return to the fixed ambient Euclidean norm at the end of
the section.

\subsection{The ambient positive definite forms}

Let $\mathcal P$ be the space of positive definite forms on $\R^N$
having determinant $\det G_p$.  Its Riemannian metric is normalized by
\[
 \mathfrak g_G(U,V)
   =\frac14\tr(G^{-1}UG^{-1}V).
\]
Here a tangent vector to $\mathcal P$ is viewed as a symmetric matrix
in any fixed ambient basis. 

\begin{lemma}\label{lem:matrixgauss}
The map $\mathcal G\colon x\mapsto G_x$ takes values in $\mathcal P$ and
satisfies
\begin{equation}\label{eq:gausspullback}
 \mathcal G^*\mathfrak g=Nh+\operatorname{Ric}_h.
\end{equation}
It is equivariant for the action of $\{g\in\mathrm{GL}(N,\R):|\det g|=1\}$ given by $g\cdot G=g^{-T}Gg^{-1}$.  Moreover, for every cone automorphism $g$ with $|\det g|=1$, we have
\begin{equation}\label{eq:gauss-cartan-distance}
 d_{\mathcal P}(G_p,G_{gp})^2
   =\sum_{i=1}^N\mu^p_i(g)^2.
\end{equation}
\end{lemma}

\begin{proof}
Use the bundle identification
\[
 \R\oplus T\Omega\longrightarrow\Omega\times\R^N,\qquad
 (\alpha,Y)\longmapsto\alpha f+f_*Y.
\]
In this identification, $G$ is the product of the standard metric on
$\R$ and $h$ on $T\Omega$.  The structure equation
\eqref{eq:structure} gives
\[
 D_X(\alpha,Y)
   =\bigl(X\alpha+h(X,Y),\,
          \nabla_XY+A_XY+\alpha X\bigr).
\]
Consequently
\[
 D=\bar\nabla+\Psi,\qquad
 \bar\nabla_X(\alpha,Y)=(X\alpha,\nabla_XY),\qquad
 \Psi_X=\begin{pmatrix}0&X^\flat\\ X&A_X\end{pmatrix}.
\]
The connection $\bar\nabla$ preserves $G$.  Cubic symmetry makes
$\Psi_X$ self-adjoint for $G$, and $\tr\Psi_X=\tr A_X=0$.  Thus, in a flat ambient trivialization,
\[
 D_XG=-2G\Psi_X,\qquad
 X(\log\det G)=-2\tr\Psi_X=0.
\]
Therefore, for every tangent vector $X$,
\begin{align*}
 \|d\mathcal G(X)\|_{\mathfrak g}^2
   &=\tr(\Psi_X^2)\\
   &=2h(X,X)+\tr(A_X^2)\\
   &=Nh(X,X)+\operatorname{Ric}_h(X,X),
\end{align*}
where the last equality is~\eqref{eq:ricci}.  Polarization proves
\eqref{eq:gausspullback}.  In particular,
\begin{equation}\label{eq:gauss-speed-bound}
 \|d\mathcal G(X)\|_{\mathfrak g}^2
   \leq(N-\varepsilon)h(X,X).
\end{equation}

The canonical affine sphere and its Blaschke metric are preserved by
every cone-preserving matrix of absolute determinant one.  The
definition of $G_x$ therefore gives
\[
 G_{gx}=g^{-T}G_xg^{-1}.
\]
This is the asserted equivariance.

We recall the distance formula with our normalization.  After a
linear change of basis taking one positive definite form to the
identity, the distance from $\id$ to a positive definite matrix $G_0$
is
\[
 d_{\mathcal P}(\id,G_0)^2
   =\frac14\sum_i(\log\lambda_i(G_0))^2.
\]
Choose a $G_p$-orthonormal ambient basis.  In this basis $G_p=\id$ and
$G_{gp}=g^{-T}g^{-1}$.  The eigenvalues of the latter matrix are the
inverse squared singular values of $g$.  The distance formula hence
gives
\[
 d_{\mathcal P}(G_p,G_{gp})^2
  =\frac14\sum_i(-2\mu_i^p(g))^2,
\]
which proves~\eqref{eq:gauss-cartan-distance}.
\end{proof}

\subsection{Complete projective chords}

We next control the length of a projective segment from $p$ to $gp$.
The first step records why completeness supplies a differential
equation on the whole real line.

\begin{lemma}\label{lem:chord-equation}
Every complete projective chord in $\Omega$ admits a Blaschke
arclength parametrization $c\colon\R\to\Omega$.  Given two distinct
points $p,q$ on that chord, choose its orientation so that $p$ precedes
$q$.  Let $[a]$ be the endpoint of the chord beyond
$q$ in the forward direction and $[b]$ its other endpoint, beyond $p$;
choose unit lifts $a,b\in\partial C\setminus\{0\}$ of these endpoints.  There are smooth positive functions $x,y$ on $\R$
such that
\[
 f(c(s))=x(s)a+y(s)b,\qquad x'(s)>0,\quad y'(s)<0.
\]
Writing $T=\dot c$, the affine-sphere lift satisfies
\begin{equation}\label{eq:chord-ode}
 \frac{d^2}{ds^2}(f\circ c)
   =\vartheta(s)\frac{d}{ds}(f\circ c)+f\circ c,\qquad
 \vartheta(s)=A(T(s),T(s),T(s)),\qquad
 |\vartheta(s)|\leq\sqrt{n-1}.
\end{equation}
\end{lemma}

\begin{proof}
In a bounded affine realization of $\Omega$, the chord is an open
interval whose endpoints lie on $\partial\Omega$.  A regular
parametrization can be reparametrized by $h$-arclength on an interval
$(s_-,s_+)$.  Suppose, for example, that $s_+<\infty$.  The curve is
then Cauchy as $s$ tends to $s_+$, since its Riemannian distance
between two points is at most their intervening arclength.
Completeness of $h$ would give a limit in $\Omega$.  This contradicts
convergence of the projective chord to its boundary endpoint.
Thus $s_+=+\infty$, and the same argument gives $s_-=-\infty$.

Every point of the chord has a unique positive expression in the
rays of $a$ and $b$, so $f\circ c=xa+yb$ for smooth $x,y>0$.  The
canonical conormal belongs to the interior dual cone.  In particular,
\[
 \nu_{c(s)}(a)>0,\qquad \nu_{c(s)}(b)>0.
\]
The conormal identities give
\[
 0=\nu_{c(s)}((f\circ c)')
   =x'\nu_{c(s)}(a)+y'\nu_{c(s)}(b).
\]
The derivative $(f\circ c)'$ is nonzero because $c$ has unit speed.
Thus $x',y'$ are nonzero and have opposite signs.  The projective
coordinate $x/y$ increases in the direction toward $[a]$, and its
derivative cannot vanish for a regular chord.  Therefore
\[
 \frac{d}{ds}(x/y)=\frac{x'y-xy'}{y^2}>0
\]
selects $x'>0$ and $y'<0$.

The vectors $f\circ c$ and $(f\circ c)'$ span the fixed ambient
two-plane of the chord.  The structure equation gives
\[
 (f\circ c)''
   =f_*(\nabla_TT+A_TT)+f\circ c.
\]
Since $f$ is transverse to the affine-sphere tangent space, the
intersection of this two-plane with that tangent space is precisely
the line spanned by $(f\circ c)'$.  It follows that
\[
 \nabla_TT+A_TT=\vartheta T.
\]
Pairing with the unit vector $T$ gives
$\vartheta=A(T,T,T)$, because
$h(\nabla_TT,T)=0$.  Finally
Proposition~\ref{prop:calabi} gives
\[
 |\vartheta|
   \leq\|A_T\|_{\mathrm{HS}}
   =\sqrt{\tr(A_T^2)}
   \leq\sqrt{n-1}.
\]
This proves~\eqref{eq:chord-ode}.
\end{proof}

The next estimate uses only the bound on the cubic coefficient.  In
particular, the strict Ricci inequality is not needed to control the
following scalar differential equations.

\begin{lemma}\label{lem:riccati}
Let $\vartheta\colon\R\to\R$ be continuous with
$|\vartheta|\leq\kappa$.  Suppose that positive $C^1$ functions $u,v$
on $\R$ satisfy
\[
 u'=1+\vartheta u-u^2,\qquad
 v'=-1+\vartheta v+v^2.
\]
Then
\[
 u_-\leq u,v\leq u_+,\qquad
 u_+=\frac{\kappa+\sqrt{\kappa^2+4}}2,\qquad
 u_-=u_+^{-1}.
\]
Consequently
\begin{equation}\label{eq:chord-ratio}
 c_\kappa\leq\frac{u}{u+v}\leq1-c_\kappa,\qquad
 c_\kappa=\frac{u_-}{u_-+u_+}>0.
\end{equation}
\end{lemma}

\begin{proof}
The constants satisfy $u_+u_-=1$ and $u_+-u_-=\kappa$, so
\[
 1+\kappa z-z^2=-(z-u_+)(z+u_-),\qquad
 1-\kappa z-z^2=(u_--z)(z+u_+).
\]
Suppose that $u(s_0)>u_+$, and reverse time by setting
$z(\tau)=u(s_0-\tau)$.  The differential inequality for $u$ gives
\[
 z'\geq(z-u_+)(z+u_-)>0
\]
as long as $z>u_+$.  In particular $z$ increases and blows up in finite time, contradicting its existence on all of $\R$. \\
If $0<u(s_0)<u_-$, the reversed function satisfies
\[
 z'\leq-(u_--z)(z+u_+).
\]
It decreases while $0<z<u_-$.  On $0<z\leq u(s_0)$ its derivative is
at most $-(u_--u(s_0))u_+<0$, so it reaches zero in finite time.
This contradicts positivity.  These two arguments prove
$u_-\leq u\leq u_+$.

The corresponding contradictions for $v$ occur in forward time.
If $v(s_0)>u_+$, then
\[
 v'\geq(v-u_+)(v+u_-)>0,
\]
which forces finite-time blow-up.  If $0<v(s_0)<u_-$, then
\[
 v'\leq-(u_--v)(v+u_+),
\]
which forces a finite-time zero.  Thus $u_-\leq v\leq u_+$ as well.
The bounds in~\eqref{eq:chord-ratio} follow by taking the minimum
and maximum of $u/(u+v)$ on $[u_-,u_+]^2$.
\end{proof}

For $p,q\in\Omega$, denote by $L(p,q)$ the Blaschke length of their
projective segment.  We emphasize that this is the length of that
particular segment, rather than the Blaschke distance.

\begin{lemma}\label{lem:chordlength}
There is a positive constant $c_p$, depending only on the fixed cone
and basepoint, such that every cone-preserving automorphism $g$ with
$|\det g|=1$ satisfies
\begin{equation}\label{eq:chordlength}
 L(p,gp)^2
 \leq\bigl(\mu_1^p(g)+c_p\bigr)
       \bigl(-\mu_N^p(g)+c_p\bigr).
\end{equation}
\end{lemma}

\begin{proof}
Assume first that $q=gp\ne p$.  Orient the complete chord from
$p$ toward $q$, so that its arclength parameter satisfies
$s(p)<s(q)$.  In Lemma~\ref{lem:chord-equation}, $[a]$ is then the
boundary endpoint beyond $q$, while $[b]$ lies beyond $p$; normalize
$a,b$ to have unit length in $\|\cdot\|_p$.  All integrals
$\int_p^q$ below are consequently taken in the positive $s$-direction.  Set
\[
 u=\frac{x'}x>0,\qquad v=-\frac{y'}y>0.
\]
Equation~\eqref{eq:chord-ode} gives the system in
Lemma~\ref{lem:riccati}, with $\kappa=\sqrt{n-1}$.  In particular,
$\zeta=u/(u+v)$ satisfies~\eqref{eq:chord-ratio}.  Direct differentiation
of the two equations gives
\[
 (\log \zeta)'
   =\frac{u'}u-\frac{u'+v'}{u+v}
   =\left(\vartheta+\frac1u-u\right)
      -(\vartheta+v-u)
   =\frac1u-v.
\]
Define
\[
 \alpha=\log\frac{x(q)}{x(p)}=\int_p^q u\,ds>0,\qquad
 \beta=\log\frac{y(q)}{y(p)}=-\int_p^q v\,ds<0.
\]
It follows that
\[
 \int_p^q\frac{ds}{u}
   =-\beta+\log\frac{\zeta(q)}{\zeta(p)}
   \leq-\beta+c_{\mathrm{ch}},\qquad
 c_{\mathrm{ch}}=\log\frac{1-c_\kappa}{c_\kappa}.
\]
Cauchy--Schwarz therefore gives
\begin{equation}\label{eq:chord-xy-length}
 L(p,q)^2
   \leq\left(\int_p^q u\,ds\right)
          \left(\int_p^q\frac{ds}{u}\right)
   \leq \alpha(-\beta+c_{\mathrm{ch}}).
\end{equation}

\noindent We next give bounds on the coefficients $\alpha$ and $\beta$ at the fixed basepoint $p$.
Choose $\ell\in C^*$, and let
\[
 \mathbb S_p=\{w\in\R^N:\|w\|_p=1\},\qquad
 m_\ell=\min_{\overline C\cap \mathbb S_p}\ell>0,\qquad
 d_p=\operatorname{dist}_{\|\cdot\|_p}(f(p),\partial C)>0.
\]
Since $\overline C$ is pointed and full-dimensional, its dual cone has
nonempty interior.  Thus a strictly positive $\ell$ can be chosen;
$\overline C\cap \mathbb S_p$ is compact, and $\ell$ has a positive minimum
there.  
Applying $\ell$ to $f(p)=x(p)a+y(p)b$ gives
\[
 x(p),y(p)\leq U_p:=\frac{\ell(f(p))}{m_\ell}.
\]
For a lower bound, the vector $y(p)b$ lies in $\partial C$, and its
distance from $f(p)$ is $x(p)$, since $\|a\|_p=1$.  Interchanging
the endpoints gives
\begin{equation}\label{eq:chord-base-coefficients}
 d_p\leq x(p),y(p)\leq U_p.
\end{equation}

\noindent At $q$, positivity gives
\[
 m_\ell x(q)\leq\ell(f(q))
   \leq\|\ell\|_{p,*}\|g f(p)\|_p
   \leq\|\ell\|_{p,*}e^{\mu_1^p(g)}.
\]
Here $\|f(p)\|_p=1$ by definition.  Combining this with
\eqref{eq:chord-base-coefficients} yields
\begin{equation}\label{eq:chord-upper-endpoint}
 \alpha\leq\mu_1^p(g)
       +\log\frac{\|\ell\|_{p,*}}{m_\ell d_p}.
\end{equation}

For the other endpoint, the two conormal equations are
\[
 x(q)\nu_q(a)+y(q)\nu_q(b)=1,\qquad
 u(q)x(q)\nu_q(a)-v(q)y(q)\nu_q(b)=0.
\]
Solving them gives
\[
 y(q)\nu_q(b)=\frac{u(q)}{u(q)+v(q)}\geq c_\kappa.
\]
The functional $\nu_p$ is the $G_p$-dual of the unit vector $f(p)$,
so $\|\nu_p\|_{p,*}=1$.  Conormal equivariance and $\|b\|_p=1$
therefore imply
\[
 \nu_q(b)\leq\|\nu_q\|_{p,*}
   =\|g^{-T}\nu_p\|_{p,*}
   \leq e^{-\mu_N^p(g)}.
\]
Consequently $y(q)\geq c_\kappa e^{\mu_N^p(g)}$.  Using
$y(p)\leq U_p$ gives
\begin{equation}\label{eq:chord-lower-endpoint}
 -\beta\leq-\mu_N^p(g)+\log(U_p/c_\kappa).
\end{equation}

The determinant normalization ensures
$\mu_1^p(g)\geq0$ and $-\mu_N^p(g)\geq0$.  We may therefore choose
one positive constant $c_p$ large enough to absorb the additive
constants in~\eqref{eq:chord-upper-endpoint} and
\eqref{eq:chord-lower-endpoint}, together with $c_{\mathrm{ch}}$ in
\eqref{eq:chord-xy-length}.  This proves~\eqref{eq:chordlength}.
If $gp=p$, its left side is zero, and the same positive choice of
$c_p$ proves the inequality directly.
\end{proof}

\subsection{The Cartan inequality}

Combining the two length comparisons proves the main estimate of the
section.

\begin{proposition}\label{prop:cartan}
There is a constant $c_p>0$ such that every cone-automorphism
$g$ normalized by $|\det g|=1$ satisfies
\begin{equation}\label{eq:cartanineq}
 \sum_{i=1}^N\mu_i^p(g)^2
 \leq(N-\varepsilon)
       \bigl(\mu_1^p(g)+c_p\bigr)
       \bigl(-\mu_N^p(g)+c_p\bigr).
\end{equation}
The constants depend only on the fixed domain, the basepoint, and
the Ricci bound, and not on $g$.
\end{proposition}

\begin{proof}
Suppose that $gp\ne p$, and apply
Lemma~\ref{lem:matrixgauss} to the projective segment from $p$ to
$gp$.  Equation~\eqref{eq:gauss-speed-bound} bounds the length of its
Gauss image by $\sqrt{N-\varepsilon}\,L(p,gp)$.  Distance is at most
the length of this image, so~\eqref{eq:gauss-cartan-distance} gives
\[
 \sum_i\mu_i^p(g)^2\leq(N-\varepsilon)L(p,gp)^2.
\]
Lemma~\ref{lem:chordlength} proves~\eqref{eq:cartanineq}.

If $gp=p$, equivariance gives $G_p=g^{-T}G_pg^{-1}$, so $g$ is
$G_p$-orthogonal.  All its logarithmic singular values vanish, and
the result follows.
\end{proof}

We finally express this estimate as a linear gap from all two-level
Cartan rays.  

\begin{proposition}\label{prop:gap}
Let $\mu(\gamma)$ be the vector of the logarithms of the singular values of the cone-automorphism $\gamma$, computed in any fixed ambient Euclidean norm, arranged in decreasing order.
There are constants $\delta>0$ and $c_{\mathrm{gap}}<\infty$ such that
\begin{equation}\label{eq:gap}
 D_k(\mu(\gamma))\geq\delta R(\mu(\gamma))-c_{\mathrm{gap}}
  \ \ \ \ \ \  \text{for all \ $\gamma\in\Gamma,\ 1\leq k<N$}.
\end{equation}
One may take $\delta=(N-1)\varepsilon/N^2$, with the additive
constant depending on the chosen norm.
\end{proposition}

\begin{proof}
We first work in the norm $G_p$.
Write
\[
 a=\mu_1^p\geq0,\qquad b=-\mu_N^p\geq0,\qquad R=a+b.
\]
If $R>0$, trace zero gives
\[
 a\leq(N-1)b,\qquad b\leq(N-1)a.
\]
Hence $a/R$ lies in $[1/N,(N-1)/N]$, and
\begin{equation}\label{eq:cartan-extreme-product}
 ab\geq\eta R^2,\qquad \eta=\frac{N-1}{N^2}.
\end{equation}
The following identity uses $\sum_i\mu^p_i=0$:
\[
 Nab-\sum_i(\mu_i^p)^2
   =\sum_i(a-\mu_i^p)(\mu_i^p+b).
\]
Every summand is nonnegative.  For $i\leq k$ it is at most
$R(a-\mu_i^p)$, whereas for $i>k$ it is at most $R(\mu_i^p+b)$.
It follows that
\begin{equation}\label{eq:cartan-defect}
 Nab-\sum_i(\mu_i^p)^2\leq R D_k(\mu^p).
\end{equation}
Expanding~\eqref{eq:cartanineq} gives constants $k_1,k_0\geq0$,
independent of $\gamma$ and $k$, such that
\[
 \sum_i(\mu_i^p)^2\leq(N-\varepsilon)ab+k_1R+k_0.
\]
Together with~\eqref{eq:cartan-extreme-product} and
\eqref{eq:cartan-defect}, this yields
\[
 R D_k(\mu^p)
   \geq\varepsilon ab-k_1R-k_0
   \geq\varepsilon\eta R^2-k_1R-k_0.
\]
For $R\geq1$, division gives
\[
 D_k(\mu^p)\geq\varepsilon\eta R-(k_1+k_0).
\]
For $0\leq R<1$, the nonnegativity of $D_k$ gives the same bound
after enlarging one additive constant.  Thus~\eqref{eq:gap} holds
in $G_p$.

To transfer the conclusion to a fixed ambient norm $\|\cdot\|_0$,
choose a linear isomorphism $S_0$ such that
$\|v\|_0=\|S_0v\|_p$.  The Cartan vector in $\|\cdot\|_0$ is the
Cartan vector of $S_0\gamma S_0^{-1}$ in $\|\cdot\|_p$.  The
singular-value inequalities for left and right multiplication give
\[
 |\mu_i^0(\gamma)-\mu_i^p(\gamma)|
   \leq e_{\mathrm{bas}},\qquad
 e_{\mathrm{bas}}=\log\bigl(\|S_0\|_p\|S_0^{-1}\|_p\bigr).
\]
In particular,
\[
 |R(\mu^0)-R(\mu^p)|\leq2e_{\mathrm{bas}},\qquad
 |D_k(\mu^0)-D_k(\mu^p)|\leq2Ne_{\mathrm{bas}}.
\]
The gap in $G_p$ therefore implies
\[
 D_k(\mu^0(\gamma))
 \geq\delta R(\mu^0(\gamma))
       -\bigl(c_{\mathrm{gap}}+2Ne_{\mathrm{bas}}+2\delta e_{\mathrm{bas}}\bigr).
\]
This proves the proposition with the same positive $\delta$.
\end{proof}

\begin{remark}\label{rem:dual-cartan-gap}
The estimate passes directly to the contragredient action.  In the
dual norm, the ordered Cartan vector of $\gamma^{-T}$ is
\[
 \widetilde\mu_i=-\mu_{N+1-i}.
\]
Consequently
\[
 R(\widetilde\mu)=R(\mu),\qquad
 D_k(\widetilde\mu)=D_{N-k}(\mu).
\]
Thus the dual group satisfies~\eqref{eq:gap} with the same constants
when the norms are dual, and with an adjusted additive constant in
any other fixed norm.
\end{remark}

\section{Projective packing in logarithmic coordinates}\label{sec:packing}

Let $n\geq2$ and $N=n+1$.  We use the Cartan vector $\mu$, its spread
$R$, and the functions $D_k$ defined above.  Matrix norms in this section are taken in the fixed standard basis.  In
particular, the singular values of a positive diagonal matrix are its
diagonal entries, arranged in decreasing order.

\begin{theorem}\label{thm:logpacking}
Let $K\subset\RP^n$ be a properly convex domain and let
$\Gamma<\operatorname{PGL}(N,\R)$ be a discrete group acting
cocompactly on $K$. Fix an affine chart $\R^n\subset\RP^n$ with positive orthant
$\R^n_{>0}$. Assume the Cartan
gap~\eqref{eq:gap}: there are constants $\delta>0$ and $c_{\mathrm{gap}}\geq0$ such that
\[
 D_k(\mu(\gamma))\geq\delta R(\mu(\gamma))-c_{\mathrm{gap}}
 \qquad
 \text{for all \ $\gamma\in\Gamma,\ 1\leq k<N$}.
\]
Then there is a constant $c_{\mathrm{pack}}$, depending only on these fixed data, such
that
\begin{equation}\label{eq:logpacking}
 \sup_{\substack{S\in\operatorname{PGL}_N(\R)\\
                  S K\subset\R_{>0}^n}}
 \int_{S K}\frac{dy_1\cdots dy_n}{y_1\cdots y_n}\leq c_{\mathrm{pack}}.
\end{equation}
\end{theorem}

We first make one fixed choice of a bounded affine realization of $K$.  This entails no restriction on the theorem.  Indeed, after
a fixed projective change $S_b$, replace $K,\Gamma,S$ by
$S_b K,S_b\Gamma S_b^{-1},S S_b^{-1}$ respectively.  The target image $S K$ is
unchanged.  Lemma~\ref{lem:cartanperturb} shows that conjugation by the
fixed matrix $S_b$ changes every Cartan coordinate by a bounded amount.
Since $R$ and $D_k$ are Lipschitz functions of these coordinates, the gap
retains the same positive $\delta$ after increasing $c_{\mathrm{gap}}$.  We therefore
work with a fixed bounded source $K\subset\R^n$, and rename the constants.
All constants below may depend on $n$, this bounded realization of $K$,
a fixed compact orbit section $\mathcal F$, and the constants in the gap, but
never on $S$.

Every image $\mathcal D=S K$ contained in the positive affine orthant is ordinarily
convex, including when it is unbounded.  For $a\in\R_{>0}^n$ and $0<r\leq1/4$, write
\[
 \mathcal B(a,r)=\prod_{i=1}^n[a_i(1-r),a_i(1+r)]
\]
for the closed relative coordinate box.  The following elementary
consequence of John's ellipsoid theorem (\cite{John48}) gives such boxes from logarithmic mass.

\begin{lemma}\label{lem:packingcellbox}
There is $c_*>0$, depending only on $n$, with the following property.
For $m\in\mathbb Z^n$, let
\[
 t_i=2^{-m_i},\qquad
 Q_m=\prod_{i=1}^n(t_i,2t_i),\qquad
 q_m=\int_{\mathcal D\cap Q_m}\frac{dy_1\cdots dy_n}{y_1\cdots y_n},
\]
where $\mathcal D\subset\R_{>0}^n$ is any convex open set.  Put
$q_{\max}=(\log2)^n<1$.  Then $q_m\leq q_{\max}$, and if
$q_m\geq\rho>0$, there is $a\in Q_m$ such that
\begin{equation}\label{eq:packingcellbox}
 \mathcal B(a,r)\subset \mathcal D,
 \qquad r=c_*\rho\leq1/4.
\end{equation}
\end{lemma}

\begin{proof}
The upper bound is the integral over the entire dyadic box:
\[
 q_m\leq\prod_{i=1}^n\int_{t_i}^{2t_i}\frac{dy_i}{y_i}
       =(\log2)^n.
\]
Assume now $q_m\geq\rho>0$. Scale $y_i=t_iY_i$, and let
$W_m=\diag(t_1^{-1},\ldots,t_n^{-1})(\mathcal D\cap Q_m)$.  It is convex and open,
is contained in $(1,2)^n$, and satisfies
\[
 \rho\leq q_m=\int_{W_m}\frac{dY_1\cdots dY_n}{Y_1\cdots Y_n}
       \leq\Vol_{\mathrm{Euc}}(W_m).
\]
Since $W_m$ is an open convex set of positive volume, its closure is a
convex body with $\operatorname{int}\overline{W_m}=W_m$ and
$\Vol_{\mathrm{Euc}}(\overline{W_m})=\Vol_{\mathrm{Euc}}(W_m)$.
The first identity is the standard fact that an open convex set is the
interior of its closure; the second follows because the boundary of a
full-dimensional convex body has Lebesgue measure zero.
Apply John's ellipsoid
theorem~\cite{John48} to the compact convex body $\overline{W_m}$.  Its
maximal inscribed ellipsoid $\mathcal E_m$, with semiaxes
$a_1',\ldots,a_n'$, satisfies
\[
 \Vol_{\mathrm{Euc}}(\mathcal E_m)
 \geq n^{-n}\Vol_{\mathrm{Euc}}(\overline{W_m})
 \geq n^{-n}\rho.
\]
Every semiaxis is at most $\sqrt n/2$, since
$\mathcal E_m\subset[1,2]^n$.  Set
\[
 c_n=\frac{n^{-n}}{\omega_n(\sqrt n/2)^{n-1}},
 \qquad
 c_*=\min\left\{\frac14,\frac{c_n}{4\sqrt n}\right\}
\]
with $\omega_n$ being the Euclidean unit-ball volume. Then, we have
\[
 \min_i a_i'
 \geq\frac{n^{-n}\rho}
 {\omega_n(\sqrt n/2)^{n-1}}
 =c_n\rho.
\]
To obtain a closed set contained in the open cell, shrink $\mathcal E_m$
by the factor $1/2$ about its centre $z$.  The resulting closed ellipsoid
lies in $\operatorname{int}\mathcal E_m\subset
\operatorname{int}\overline{W_m}=W_m$.  It contains a closed Euclidean ball
of radius $c_n\rho/2$, and hence a closed coordinate cube about $z$ with
half-width $c_n\rho/(2\sqrt n)$.

Set $a_i=t_i z_i$, so $t_i<a_i<2t_i$.  Scaling the cube back gives
coordinate half-widths $c_n\rho t_i/(2\sqrt n)$.  With $c_*$ fixed as above, set $r=c_*\rho$.
Then $r a_i\leq2c_*\rho t_i\leq c_n\rho t_i/(2\sqrt n)$, which proves
the box inclusion.  Since $\rho\leq q_{\max}<1$, we also have $r\leq1/4$.
\end{proof}

\begin{lemma}[Relative-box normalization]\label{lem:packingnormalization}
Suppose $\mathcal D=S K\subset\R_{>0}^n$ contains
$\mathcal B(a,r)$ and $\mathcal B(b,r)$, where $0<r\leq1/4$.  For
$0\leq\lambda\leq1$ put
\[
 p^\lambda=(1-\lambda)a+\lambda b,
 \qquad
 \mathsf D_\lambda=\diag\bigl((p_1^\lambda)^{-1},\ldots,
                      (p_n^\lambda)^{-1},1\bigr),
\]
and define
\begin{equation}\label{eq:packingnormalizers}
 \mathcal J(Y)=\frac{Y}{1+\sum_iY_i},
 \qquad \Phi_\lambda=\mathcal J \mathsf D_\lambda S.
\end{equation}
Let
\[
 \mathcal T=\{z_i>0,\ \textstyle\sum_i z_i<1\},
 \qquad b_*=\frac1N(1,\ldots,1),
 \qquad \beta_n=\frac1{8N(1+\sqrt n)}.
\]
Then for every $\lambda$,
\begin{equation}\label{eq:packinginnerouter}
 \overline B_{\mathrm{Euc}}(b_*,\beta_n r)
 \subset \Phi_\lambda K\subset\mathcal T.
\end{equation}
\end{lemma}

\begin{proof}
Minkowski addition of coordinate intervals gives the exact identity
\[
 (1-\lambda)\mathcal B(a,r)+\lambda\mathcal B(b,r)
       =\mathcal B(p^\lambda,r).
\]
Thus convexity of $\mathcal D$ implies
$\mathcal B(p^\lambda,r)\subset \mathcal D$.  After applying $\mathsf D_\lambda$, this is
the same box $[1-r,1+r]^n$ for every parameter.  The entire positive
orthant is mapped by $\mathcal J$ into $\mathcal T$, which proves the outer
inclusion even if $\mathcal D$ is unbounded.

For the inner inclusion use the explicit inverse
\[
 \mathcal J^{-1}(z)=\frac{z}{1-\sum_i z_i}.
\]
Writing $\boldsymbol1=(1,\ldots,1)$ and $z=b_*+w$, we obtain
\[
 \mathcal J^{-1}(b_*+w)-\boldsymbol1
   =\frac{w+\boldsymbol1\sum_i w_i}{1/N-\sum_i w_i}.
\]
If $\|w\|_2\leq\beta_n r$ and $r\leq1/4$, then
$w_i>-1/N$ and
$1-\sum_i(b_{*,i}+w_i)=1/N-\sum_i w_i>0$.
Hence $b_*+w\in\mathcal T$, and $\mathcal J^{-1}$ is defined there.
Moreover the denominator in the displayed formula is at least
$1/(2N)$, while the numerator in the sup norm is at most
$(1+\sqrt n)\|w\|_2$.  Therefore
\[
 \|\mathcal J^{-1}(b_*+w)-\boldsymbol1\|_\infty
 \leq 2N(1+\sqrt n)\beta_n r=r/4<r.
\]
The asserted ball
is consequently contained in the image of $[1-r,1+r]^n$, and hence in
$\Phi_\lambda K$.
\end{proof}

\begin{lemma}\label{lem:packingcorrection}
There is a constant $c_{\mathrm{corr}}\geq1$, depending only on
$n$, the bounded realization of $K$, a compact fundamental domain $\mathcal F$ and
$R_V=\operatorname{diam}\mathcal T$ such that the following holds.  If $0<r\leq1/4$ and
a projective map $\Phi'$ satisfies
\[
 \overline B_{\mathrm{Euc}}(b_*,\beta_n r)
 \subset \Phi'K\subset\mathcal T,
\]
one can choose $\gamma'\in\Gamma$ so that the projective map
$\mathsf W'=\Phi'(\gamma')^{-1}$ has condition number
\begin{equation}\label{eq:packingcorrectioncondition}
 \operatorname{cond}(\mathsf W'):=\|\mathsf W'\| \|(\mathsf W')^{-1}\|\leq c_{\mathrm{corr}}/r.
\end{equation}
Its absolute-determinant-one representative and its inverse both have
operator norm at most $c_{\mathrm{corr}}/r$.
Moreover, if $\Phi',\Phi''$ satisfy the same hypothesis, then
\begin{equation}\label{eq:packingcartancomparison}
 \bigl\|\mu(\Phi'(\Phi'')^{-1})
       -\mu(\gamma'(\gamma'')^{-1})\bigr\|_\infty
 \leq2\log(c_{\mathrm{corr}}/r).
\end{equation}
\end{lemma}

\begin{proof}
Fix a compact fundamental domain $\mathcal F\subset \mathrm{int}(K)$. Since
$b_*\in\Phi'K$, the point $(\Phi')^{-1}b_*$ belongs to $K$.
Cocompactness therefore supplies $\gamma'\in\Gamma$ such that
$x=\gamma'(\Phi')^{-1}b_*\in \mathcal F$.  Then $\mathsf W'$ maps $K$ to $\Phi'K$ and
sends $x$ to $b_*$.  There are fixed positive constants $r_U,R_U,R_V$
such that the source and target translates
\[
 \mathcal U=K-x,\qquad \mathcal V=\Phi'K-b_*
\]
satisfy
\[
 \begin{gathered}
 B_{\mathrm{Euc}}(0,r_U)\subset \mathcal U\subset B_{\mathrm{Euc}}(0,R_U),
 \\
 B_{\mathrm{Euc}}(0,\beta_n r)\subset \mathcal V
       \subset B_{\mathrm{Euc}}(0,R_V).
 \end{gathered}
\]
We may take $R_V=\operatorname{diam}\mathcal T$ after translating by $b_*$.

We compute condition numbers using any invertible matrix representative:
scalar multiplication does not change them. Translate $x$ and $b_*$
to zero. An invertible matrix representing the
resulting projective map has blocks
$\left(\begin{smallmatrix}\mathsf L&a\\c^T&d\end{smallmatrix}\right)$.
Since zero maps to the finite affine point zero, $a=0$ and $d\ne0$.
After dividing the whole matrix by $d$ and renaming its blocks,
we obtain the representative normalized by its bottom-right entry $d=1$:
\[
 \widehat{\mathsf W}=\begin{pmatrix}\mathsf L&0\\c^T&1\end{pmatrix},
 \qquad u\longmapsto\frac{\mathsf L u}{1+c\cdot u}.
\]
The denominator is positive on $\mathcal U$: it equals one at zero, it has no
zero in $\mathcal U$ because the image is entirely finite, and $\mathcal U$ is connected.
Its positivity on $B_{\mathrm{Euc}}(0,r_U)$ implies $\|c\|\leq r_U^{-1}$.  For a
unit vector $v$, evaluate at $u=(r_U/2)v$.  The denominator belongs to
$[1/2,3/2]$ and the image has norm at most $R_V$, so
\[
 \|\mathsf L v\|\leq3R_V/r_U.
\]
It follows, for this particular representative, that
\begin{equation}\label{eq:packingforwardmatrix}
 \|\widehat{\mathsf W}\|\leq1+(1+3R_V)/r_U\leq c_{\mathrm{fwd}}.
\end{equation}
The inverse is
\[
 \widehat{\mathsf W}^{-1}
   =\begin{pmatrix}\mathsf L^{-1}&0\\-c^T\mathsf L^{-1}&1\end{pmatrix}.
\]
Its own denominator is positive on $\mathcal V$, it fixes zero, and its image is
contained in $\mathcal U$.  Applying positivity on its source ball of radius
$\beta_n r$ gives
$\|c^T\mathsf L^{-1}\|\leq(\beta_n r)^{-1}$.  Evaluation at
$v\beta_n r/2$, for a unit vector $v$, gives
$\|\mathsf L^{-1}v\|\leq3R_U/(\beta_n r)$.  Consequently
\begin{equation}\label{eq:packinginversematrix}
 \|\widehat{\mathsf W}^{-1}\|
 \leq1+(1+3R_U)/(\beta_n r)\leq c_{\mathrm{inv}}/r.
\end{equation}
The two translations and their inverses have uniformly bounded
matrix representatives, because $x\in \mathcal F$ and $b_*$ is fixed.  Undoing these
translations gives a representative $\mathsf W$ of $\mathsf W'$ with
$\|\mathsf W\|\leq c_{\mathrm{fwd}}'$ and
$\|\mathsf W^{-1}\|\leq c_{\mathrm{inv}}'/r$.  In particular,
\[
 \operatorname{cond}(\mathsf W)
 =\|\mathsf W\|\,\|\mathsf W^{-1}\|
 \leq\frac{c_{\mathrm{fwd}}' c_{\mathrm{inv}}'}{r}.
\]
This estimate is projectively well defined, since condition number is
unchanged by scalar rescaling.  Normalize $\mathsf W$ by
$\widetilde{\mathsf W}=|\det\mathsf W|^{-1/N}\mathsf W$.  The singular values of
$\widetilde{\mathsf W}$ have geometric mean one, so
$s_1(\widetilde{\mathsf W})\geq1\geq s_N(\widetilde{\mathsf W})$.  Hence
\[
 \|\widetilde{\mathsf W}\|=s_1(\widetilde{\mathsf W})
 \leq\operatorname{cond}(\mathsf W),\qquad
 \|\widetilde{\mathsf W}^{-1}\|=s_N(\widetilde{\mathsf W})^{-1}
 \leq\operatorname{cond}(\mathsf W).
\]
The existence of the constant $c_{\mathrm{corr}}\geq1$ as in the statement follows.

Finally choose absolute-determinant-one representatives of $\Phi',\Phi''$ and
$\gamma',\gamma''$, and define $\mathsf W'=\Phi'(\gamma')^{-1}$ and
$\mathsf W''=\Phi''(\gamma'')^{-1}$ as matrices.  All these matrices have absolute
determinant one, and
\[
 \Phi'(\Phi'')^{-1}
   =\mathsf W'\bigl(\gamma'(\gamma'')^{-1}\bigr)(\mathsf W'')^{-1}.
\]
The singular-value inequalities used in
Lemma~\ref{lem:cartanperturb} state that
\[
 \frac{s_i(g_0)}{\|\mathsf U^{-1}\|\,\|\mathsf V^{-1}\|}
 \leq s_i(\mathsf U g_0\mathsf V)\leq\|\mathsf U\|\,\|\mathsf V\|\,s_i(g_0).
\]
Use $\mathsf U=\mathsf W'$, $\mathsf V=(\mathsf W'')^{-1}$ and the normalized matrix bounds.
This gives \eqref{eq:packingcartancomparison}.
\end{proof}

\begin{lemma}\label{lem:packingseparation}
There is a constant $c_{\mathrm{sep}}$ such that, whenever a projective image
$\mathcal D=S K\subset\R_{>0}^n$ contains $\mathcal B(a,r)$ and
$\mathcal B(b,r)$ with $0<r\leq1/4$, one has
\begin{equation}\label{eq:packingseparation}
 \max_i|\log(b_i/a_i)|\leq c_{\mathrm{sep}}(1+\log(1/r)).
\end{equation}
\end{lemma}

\begin{proof}
Use the maps $\Phi_\lambda$ from
Lemma~\ref{lem:packingnormalization}. For each parameter, apply
Lemma~\ref{lem:packingcorrection} with the same value of $r$ to
$\Phi_\lambda$, obtaining an element $\gamma_\lambda$.
The relative Cartan comparison~\eqref{eq:packingcartancomparison}
then holds uniformly for any pair of parameters.  
Put $\xi_i=\log(b_i/a_i)$ and set $\lambda=e^{-s}$ for $s\geq\log2$.
Write $\Phi_s=\Phi_{e^{-s}}$ and $\gamma_s=\gamma_{e^{-s}}$.  Then
\[
 p_i(s)=a_i\bigl(1-e^{-s}+e^{\xi_i-s}\bigr).
\]
Since $1-e^{-s}\in[1/2,1)$, the elementary bounds for the sum of the two
positive terms give
\begin{equation}\label{eq:packingtropicalapproximation}
 \bigl|\log p_i(s)-\log a_i-(\xi_i-s)_+\bigr|\leq\log2.
\end{equation}
We first control the ratios $b_i/a_i$ that are large, namely
$\xi_i=\log(b_i/a_i)>\log2$; the opposite ratios are handled
at the end by interchanging $a$ and $b$.
Suppose that $\xi_{\max}=\max_i\xi_i>\log2$.  Subdivide the
interval $[\log2,\xi_{\max}]$ at $\xi_i$.  There are at most $n$ positive-length pieces.  On a piece
$[s_1,s_2]$ let $I=\{i:\xi_i>s\}$ for an interior point $s$, and put
$\tau=s_2-s_1$. The set $I$ is constant and nonempty on the interior.
On this piece every coordinate satisfies the dichotomy
$\xi_i\leq s_1$ or $\xi_i\geq s_2$; the second alternative
holds exactly when $i\in I$. Consequently,
\[
 (\xi_i-s_1)_+-(\xi_i-s_2)_+=\tau\mathbf1_I(i)
\]
holds also at the endpoints.  Applying
\eqref{eq:packingtropicalapproximation} there gives
\begin{equation}\label{eq:packingtropicalratios}
 \ell_i:=\log\frac{p_i(s_1)}{p_i(s_2)}
       =\tau\mathbf1_I(i)+e_i,\qquad |e_i|\leq2\log2.
\end{equation}
Write $j=|I|$, so $1\leq j\leq n<N$.
By definition of $\Phi_{s_2}$, we have
\begin{equation}\label{eq:packingrelativediagonal}
 \Phi_{s_2}\Phi_{s_1}^{-1}
 =\mathsf J\diag\left(
       \frac{p_1(s_1)}{p_1(s_2)},\ldots,
       \frac{p_n(s_1)}{p_n(s_2)},1\right)\mathsf J^{-1}.
\end{equation}
with
\[
 \mathsf J=\begin{pmatrix}I_n&0\\\boldsymbol1^T&1\end{pmatrix}.
\]
Append $\ell_N=e_N=0$ in \eqref{eq:packingtropicalratios}.  Divide the
diagonal matrix by the $N$-th root of its determinant and let
\[
 Z_j=\bigl(
   \underbrace{(N-j)/N,\ldots,(N-j)/N}_{j\text{ entries}},
   \underbrace{-j/N,\ldots,-j/N}_{N-j\text{ entries}}\bigr).
\]
Then the Cartan vector of the normalized diagonal matrix is within
$4\log2$ of $\tau Z_j$.  To bound the Cartan vector of the conjugated matrix, we let
$\mathsf F=\operatorname{diag}(1,\ldots,1,-1)$. The matrix
$\mathsf J$ has determinant one and satisfies
$\mathsf J^{-1}=\mathsf F\mathsf J\mathsf F$.
Since $\mathsf F$ is orthogonal,
$\|\mathsf J^{-1}\|=\|\mathsf J\|$.
For any absolute-determinant-one matrix $g_0$, the min--max singular-value
inequalities, applied to $\mathsf U=\mathsf J$ and $\mathsf V=\mathsf J^{-1}$, give
\[
 \left|\mu_i(\mathsf Jg_0\mathsf J^{-1})-\mu_i(g_0)\right|
 \leq\log\bigl(\|\mathsf J\|\cdot\|\mathsf J^{-1}\|\bigr)
 =\log\operatorname{cond}\mathsf J.
\]
Using
\eqref{eq:packingcartancomparison}, we conclude that
\begin{equation}\label{eq:packingtropicalerror}
 \|\mu(\gamma_{s_2}\gamma_{s_1}^{-1})-\tau Z_j\|_\infty
 \leq E(r):=4\log2+\log\operatorname{cond}\mathsf J+2\log(c_{\mathrm{corr}}/r).
\end{equation}
The vector $Z_j$ satisfies $R(Z_j)=1$ and $D_j(Z_j)=0$.  If an
ordered vector $v$ is within $E(r)$ of $\tau Z_j$, then $R(v)\geq \tau-2E(r)$.
Every summand in $D_j(v)$ is at most $2E(r)$, since
the corresponding summand for $\tau Z_j$ is zero; hence
$D_j(v)\leq2NE(r)$.  Apply the Cartan gap~\eqref{eq:gap} to
the group element $\gamma_{s_2}\gamma_{s_1}^{-1}$.  It gives
\[
 \delta\bigl(\tau-2E(r)\bigr)-c_{\mathrm{gap}}\leq2NE(r),
\]
and thus
\begin{equation}\label{eq:packingtropicalinterval}
 \tau\leq\left(2+\frac{2N}{\delta}\right)E(r)+\frac{c_{\mathrm{gap}}}{\delta}
       =:L_{\mathrm{piece}}(r).
\end{equation}
This is uniform over all active sets and all pieces.  Summing at most
$n$ pieces yields
\[
 \max_i\xi_i\leq\log2+nL_{\mathrm{piece}}(r).
\]
The same inequality holds directly when $\xi_{\max}\leq\log2$.
Interchanging $a$ and $b$ and repeating the construction bounds
$\max_i(-\xi_i)$ by the same expression.  Consequently
\[
 \max_i|\xi_i|\leq\log2+nL_{\mathrm{piece}}(r)
       \leq c_{\mathrm{sep}}(1+\log(1/r)),
\]
as required.  Breakpoints with equal values create no positive-length
piece, and the endpoint identities used above remain valid there.
\end{proof}

\begin{lemma}\label{lem:packinglayercake}
Let $(q_m)_{m\in\mathbb Z^n}$ be nonnegative numbers with
$q_m\leq q_{\max}=(\log2)^n$.  If
\[
 \#\{m:q_m\geq\rho\}\leq c_{\mathrm{cnt}}(1+|\log\rho|)^n
 \qquad( \text{with \ $0<\rho\leq q_{\max}$}),
\]
then
\[
 \sum_{m\in\mathbb Z^n}q_m
 \leq c_{\mathrm{cnt}}\int_0^{q_{\max}}(1+|\log\rho|)^n\,d\rho<\infty.
\]
\end{lemma}

\begin{proof}
Tonelli's theorem for nonnegative functions gives
\[
 \sum_{m\in \mathbb Z^n}q_m
 =\sum_m\int_0^{q_{\max}}\mathbf1_{\{\rho<q_m\}}\,d\rho
 =\int_0^{q_{\max}}\#\{m:q_m>\rho\}\,d\rho \leq c_{\mathrm{cnt}}\int_0^{q_{\max}}(1+|\log\rho|)^n\,d\rho.
\]
by assumption. Since $q_{\max}<1$, we conclude
\[
 \int_0^{q_{\max}}(1+|\log\rho|)^n\,d\rho
 \leq\int_0^1(1-\log\rho)^n\,d\rho
 =\int_0^\infty(1+t)^ne^{-t}\,dt<\infty.
\]
\end{proof}

\begin{proof}[Proof of Theorem~\ref{thm:logpacking}]
Fix $S$ and put $\mathcal D=S K\subset\R_{>0}^n$.  Use the dyadic boxes $Q_m$; they partition $\mathcal D$ up to subsets of measure zero.
Fix $0<\rho\leq q_{\max}$.  For every index with $q_m\geq\rho$, choose a
centre $a(m)\in Q_m$ such that
$\mathcal B(a(m),c_*\rho)\subset \mathcal D$.  If there are two such indices $m,m'$, then
Lemma~\ref{lem:packingseparation} gives
\[
 \max_i|\log a_i(m)-\log a_i(m')|
 \leq c_{\mathrm{sep}}\bigl(1+\log(1/(c_*\rho))\bigr).
\]
Write $a_i(m)=2^{-m_i}z_i(m)$, where $1<z_i(m)<2$.  It follows that
\begin{align*}
 |m_i-m'_i|\log2
 &\leq |\log a_i(m)-\log a_i(m')|
       +|\log z_i(m)-\log z_i(m')|\\
 &\leq c_{\mathrm{sep}}\bigl(1+\log(1/(c_*\rho))\bigr)+\log2.
\end{align*}
Because $c_*$ is fixed and $0<\rho\leq q_{\max}<1$, all high-mass indices
have pairwise lattice distance at most $c_{\mathrm{lat}}(1+|\log\rho|)$ in the sup
norm.  If the set is nonempty, fix one of its indices $m^0$; all the
others lie in an integer-coordinate cube about $m^0$ of this radius.
Consequently
\begin{equation}\label{eq:packingcellcount}
 \#\{m:q_m\geq\rho\}
 \leq\bigl(2\lfloor c_{\mathrm{lat}}(1+|\log\rho|)\rfloor+1\bigr)^n
 \leq c_{\mathrm{cnt}}(1+|\log\rho|)^n.
\end{equation}
If the set of high-mass indices is empty, the same bound is clearly satisfied.  Finally, countable additivity on the dyadic cells and
Lemma~\ref{lem:packinglayercake} give
\[
 \int_{\mathcal D}\frac{dy_1\cdots dy_n}{y_1\cdots y_n}
 =\sum_{m\in\mathbb Z^n}q_m
 \leq c_{\mathrm{cnt}}\int_0^{q_{\max}}(1+|\log\rho|)^n\,d\rho<\infty.
\]
All constants are independent of $S$, which proves the asserted
supremum.
\end{proof}

\section{Uniform volumes of projective simplices}\label{sec:volumes}

We now convert the estimate in the previous section into a bound for projective simplices.  Recall that \(\omega_n\) denotes the Euclidean
volume of the unit ball in \(\R^n\).
The Holmes--Thompson normalization is the one in
\eqref{eq:htdensity}: the density is the cotangent unit-ball volume divided by
\(\omega_n\).  We write \(\Vol_{F,L}(E)\) for the volume of \(E\) computed
with the Funk metric of the enclosing convex domain \(\operatorname{int}L\).

First, Theorem \ref{thm:logpacking} applies to the dual domain. Indeed, the canonical
conormal map induces a \(\Gamma\)-equivariant diffeomorphism
\(\Omega\longrightarrow\Omega^*\), where the action on \(\Omega^*\) is
contragredient; see \cite[Theorem~3 and Proposition~8]{Loftin08}.
Consequently this action is cocompact.  In the dual ambient norm, its Cartan
vector is
\[
  \widetilde\mu_i(\gamma)=-\mu_{N+1-i}(\gamma),
  \qquad 1\le i\le N.
\]
For the spread and the defects used in the Cartan gap this gives
\[
  R(\widetilde\mu)=R(\mu),
  \qquad D_k(\widetilde\mu)=D_{N-k}(\mu).
\]
Thus the same gap holds for the contragredient action.  We may therefore
apply Theorem~\ref{thm:logpacking} to every projective image of
\(\Omega^*\) contained in a positive orthant.\\

We will use the following precise form of Funk volume duality.  For a properly
convex body \(K\subset\RP^n\), let \(K^\vee\) be its projective dual body,
the closure of the set of hyperplanes disjoint from \(K\).

\begin{proposition}\cite[Corollary~3.11]{Faifman21}\label{prop:duality}
  Let \(K,L\subset\RP^n\) be properly convex bodies with nonempty interiors,
  and suppose \(K\subset\operatorname{int}L\).  Then
  \begin{equation}\label{eq:funkduality}
    \Vol_{F,L}(K)=\Vol_{F,K^\vee}(L^\vee).
  \end{equation}
  These volumes are independent of the affine charts used to define the
  Funk metrics.  
\end{proposition}

The density on the enclosing dual simplex can be computed directly.

\begin{lemma}\label{lem:simplexdensity}
  Let
  \[
    \mathcal S=\{(a_0,\ldots,a_n):a_i>0,\ \textstyle\sum_{i=0}^n a_i=1\}.
  \]
  In the coordinates \(a_1,\ldots,a_n\), with
  \(a_0=1-\sum_{i=1}^n a_i\), the Funk--Holmes--Thompson density of
  \(\mathcal S\) is
  \begin{equation}\label{eq:simplexdensity}
    d\Vol_{F,\mathcal S}
    =\frac{1}{n!\omega_n}
      \frac{da_1\cdots da_n}{a_0a_1\cdots a_n}.
  \end{equation}
\end{lemma}

\begin{proof}
  At \(a\in\mathcal S\), the tangent unit ball is \(\mathcal S-a\), using Minkowski-sum notation.
  Its polar, the cotangent unit ball, is the simplex with vertices
  \[
    \varpi_i=-\frac{da_i}{a_i}\quad(1\le i\le n),
    \qquad \varpi_0=\frac{da_1+\cdots+da_n}{a_0}.
  \]
  In the cotangent basis $da_1,\ldots,da_n$, the volume of this simplex is $|\det W_a|/n!$, where $-W_a$ is the matrix whose columns are $\varpi_i-\varpi_0$. Hence,
  \[
    W_a=\diag(a_1^{-1},\ldots,a_n^{-1})
       +a_0^{-1}\boldsymbol 1\boldsymbol 1^T
  \]
  and
  \[
    |\det W_a|
    =\frac{1}{a_1\cdots a_n}
       \left(1+\frac{a_1+\cdots+a_n}{a_0}\right)
    =\frac{1}{a_0\cdots a_n}.
  \]
  Dividing by \(\omega_n\) proves the formula.
\end{proof}

\begin{proof}[Proof of Theorem~\ref{thm:main}]
  Let \(\sigma\subset\Omega\) be a nondegenerate projective
  \(n\)-simplex with interior vertices.  Choose positive independent cone
  lifts \(v_0,\ldots,v_n\in C\) and put \(v_\sigma=\sum_i v_i\).
  On the dual affine section \(\ell(v_\sigma)=1\), set
  \[
    a_i=\ell(v_i),\qquad 0\le i\le n.
  \]
  These coordinates identify the dual domain with a convex domain
  \(\mathcal D_\sigma\subset\mathcal S\).  The enclosing simplex
  \(\mathcal S\) is the interior of \(\sigma^\vee\).
  Since \(\sigma\) is compactly contained in \(\Omega\), the dual
  body \(\overline{\mathcal D_\sigma}\) is compactly contained in \(\mathcal S\).
  Proposition~\ref{prop:duality} and
  Lemma~\ref{lem:simplexdensity} give the exact identity
  \begin{equation}\label{eq:dualintegral}
    \Vol_{F,\Omega}(\sigma)
    =\Vol_{F,\mathcal S}(\mathcal D_\sigma)
    =\frac{1}{n!\omega_n}
       \int_{\mathcal D_\sigma}
          \frac{da_1\cdots da_n}{a_0a_1\cdots a_n}.
  \end{equation}
  Partition \(\mathcal D_\sigma\) into Borel sets \(E_j\), \(0\le j\le n\),
  on which \(a_j\) is a largest coordinate, resolving ties by the smallest
  index.  On \(E_j\) we have \(a_j\ge1/N\), so
  \[
    \int_{E_j}\frac{da_1\cdots da_n}{\prod_{i=0}^n a_i}
    \le N\int_{\mathcal D_\sigma}
       \frac{\prod_{i\ne j}da_i}{\prod_{i\ne j}a_i}.
  \] 
  Theorem~\ref{thm:logpacking} 
  bounds the integrals by one constant \(c_{\mathrm{dual}}\), independent of
  \(\sigma\) and \(j\).  Summing yields
  \[
    \int_{\mathcal D_\sigma}\frac{da_1\cdots da_n}{\prod_{i=0}^n a_i}
    \le N\sum_{j=0}^n
       \int_{\mathcal D_\sigma}
          \frac{\prod_{i\ne j}da_i}{\prod_{i\ne j}a_i}
    \le N^2c_{\mathrm{dual}}.
  \]
  The Funk volume bound is therefore
  \[
    \Vol_{F,\Omega}(\sigma)\le\frac{N^2c_{\mathrm{dual}}}{n!\omega_n}.
  \]
  To pass to the Hilbert metric, we apply the Rogers--Shephard inequality \cite[Theorem 1]{RS57}
  \[
    d\Vol_{\mathrm{HT},\Omega}
    \le 2^{-n}\binom{2n}{n}\,d\Vol_{F,\Omega}.
  \]
   Hence
  \[
    \Vol_{\mathrm{HT},\Omega}(\sigma)
    \le 2^{-n}\binom{2n}{n}\frac{N^2c_{\mathrm{dual}}}{n!\omega_n}.
  \]
  The bound is uniform over ordinary projective simplices, as required.
  A degenerate \(n\)-simplex has dimension less than \(n\) and hence zero
  \(n\)-dimensional volume, so the same bound holds in that case.
\end{proof}

\section{Straightening and simplicial volume}\label{sec:cocycle}

The uniform bound established above is an estimate on simplices
in the convex domain $\Omega$.  A projected simplex is only a
singular simplex and may fail to be embedded; however, neither embeddedness
nor non-overlap is needed.  We use the canonical affine sphere to define an equivariant straightening homotopic to the identity, from which we will conclude the proof of Corollaries \ref{cor:positive} and \ref{cor:alternative}. The final argument is similar to the proof of the lower bound of the simplicial volume for closed hyperbolic manifolds in \cite{BenedettiPetronio}. We, therefore, only explain the main ideas and refer the interested reader to Benedetti--Petronio's book for details. \\

First suppose that \(M=\Gamma\backslash\Omega\) is oriented, and let
\(\omega_h\) be its positive smooth Blaschke volume form.  
The ratio of the Blaschke density and the Holmes--Thompson density is continuous and \(\Gamma\)-invariant.  It descends to the compact quotient,
so there is \(b_\Omega<\infty\) such that
\[
  d\Vol_h\le b_\Omega\,d\Vol_{\mathrm{HT},\Omega}
  \quad\text{on }\Omega.
\]
By Theorem~\ref{thm:main}, there is consequently \(V_h<\infty\) bounding
the Blaschke volumes of all ordinary projective \(n\)-simplices in
\(\Omega\) with interior vertices.

 
\begin{proposition}\label{prop:straightening}
  For every integer \(0\leq k\leq n+1\) and points $(x_i)_{i=0, \dots, k} \in \Omega$, the formula
  \begin{equation}\label{eq:straight}
    \st(x_0,\ldots,x_k)(t_0,\ldots,t_k)
      =\left[\sum_{i=0}^k t_i f(x_i)\right]
  \end{equation}
  defines a smooth simplex in \(\Omega\).  It is \(\Gamma\)-equivariant
  and compatible with faces.  It induces a chain map on the continuous
  singular chains of \(M\), chain-homotopic to the identity.  If the input
  singular simplex is smooth, all simplices in the homotopy may be
  taken smooth. For \(k=n\) and projectively independent vertices, the straightened simplex is the ordinary projective convex hull of its vertices.  For projectively
  dependent vertices its \(n\)-volume vanishes.
\end{proposition}

\begin{proof}
  The sum in \eqref{eq:straight} belongs to the open convex cone \(C\).
  Equivariance follows from the equivariance of \(f\), and setting one coordinate $t_{i}$ equal to zero gives the formula for the corresponding face.

  Let \(\pi:\Omega\to M\) be the covering map.  For a continuous singular
  simplex \(\sigma:\Delta^k\to M\), choose a lift
  \(\widetilde\sigma\), set \(x_i=\widetilde\sigma(e_i)\), and define
  \[
    \st_k\sigma=\pi\circ\st(x_0,\ldots,x_k).
  \]
  A different lift differs by a deck transformation, so equivariance makes
  this definition independent of the lift.  Face compatibility proves that
  it is a chain map. On each lifted simplex the map
  \[
    \mathscr H_{\tilde\sigma}(t,u)
      =\left[(1-u)f(\widetilde\sigma(t))
           +u\sum_{j=0}^k t_j f(x_j)\right],
    \qquad (t,u)\in\Delta^k\times[0,1].
  \]
  is an equivariant homotopy from
  \(\widetilde\sigma\) to \(\st(x_0,\ldots,x_k)\).  If \(\sigma\) is
  smooth, then so is this homotopy.
\end{proof}

\begin{proposition}\label{prop:cocycle}
  On continuous real singular \(n\)-chains of \(M\), the formula
  \[
    c_{\mathrm{vol}}(\sigma)=\int_{\st_n\sigma}\omega_h
  \]
  defines a bounded cocycle with \(\|c_{\mathrm{vol}}\|_\infty\le V_h\).  Its evaluation
  on the fundamental class is
  \[
    c_{\mathrm{vol}}([M])=\Vol_h(M)>0.
  \]
\end{proposition}

\begin{proof}
  The integral is defined for every continuous singular simplex, because
  its straightening is smooth. Compatibility with faces and Stokes' theorem show that
$c_{\mathrm{vol}}$ is a cocycle.
Proposition~\ref{prop:straightening} and
  the uniform bound \(V_h\) show \(|c_{\mathrm{vol}}(\sigma)|\le V_h\), including
  degenerate cases, hence $c_{\mathrm{vol}}$ is a bounded cocycle.  Choose an oriented smooth triangulation of \(M\), and let \(z_0\) be its smooth singular fundamental cycle.  It satisfies
  \(\int_{z_0}\omega_h=\int_M\omega_h=\Vol_h(M)\).
  Every real continuous singular fundamental cycle \(z\) is
  homologous to \(z_0\).  Their difference is a singular boundary, on
  which the cocycle \(c_{\mathrm{vol}}\) vanishes.  Hence
  \(c_{\mathrm{vol}}([M])=c_{\mathrm{vol}}(z_0)=\Vol_h(M)\).  
\end{proof}

\begin{proof}[Proof of Corollary~\ref{cor:positive}]
  If $M$ is oriented, for every real singular fundamental cycle \(z\),
  Proposition~\ref{prop:cocycle} gives
  \[
    \Vol_h(M)=|c_{\mathrm{vol}}(z)|\le\|c_{\mathrm{vol}}\|_\infty\,\|z\|_1.
  \]
  Taking the infimum over all such cycles yields the usual bounded-cocycle
  lower bound for simplicial volume \cite{Gromov82}:
  \begin{equation}\label{eq:lowerbound}
    \|M\|\ge\frac{\Vol_h(M)}{\|c_{\mathrm{vol}}\|_\infty}>0.
  \end{equation}

  If the original manifold is nonorientable, apply the argument to its
  orientation double cover \(\widehat M\), which has the same properly
  convex universal cover and irreducible cone.  We use the standard
  convention for nonorientable simplicial volume, under which
  \(\|M\|=\tfrac12\|\widehat M\|\).  Thus the positivity conclusion holds
  in this case as well.
\end{proof}

\begin{proof}[Proof of Corollary \ref{cor:alternative}]
We only need to prove that if the cone is reducible then $\|M\|=0$. The decomposition of a proper convex cone into indecomposable direct factors is unique up to permutation; see \cite[Chapter~II, Remark~1.11]{AMRT10}.  Consequently, cone automorphisms permute these factors.  After passing to a finite-index subgroup $\Gamma_0<\Gamma$, we may therefore assume that $\Gamma_0$ preserves every factor and the orientation of
$\Omega$. We obtain a
$\Gamma_0$-invariant decomposition
\[
 \R^N=L_1\oplus L_2,
 \qquad C=C_1\times C_2,
\]
where each $C_i\subset L_i$ is a proper open convex cone.
Set $d_i=\dim L_i>0$, and consider
\[
 g_t=e^{d_2t}\id_{L_1}\oplus e^{-d_1t}\id_{L_2},
 \qquad t\in\R.
\]
These transformations preserve $C$, have determinant one,
and commute with $\Gamma_0$.  By the naturality of the
canonical affine sphere, they preserve its Blaschke metric $h$.
Their projective action therefore descends to an isometric
flow on the closed manifold
$\widehat M=\Gamma_0\backslash\Omega$.

This flow has a nowhere vanishing generating vector field.
Indeed, every vector in $C$ has nonzero components in both
$L_1$ and $L_2$, and the two components are rescaled with
different weights.  Thus no point of $\Omega$ is stationary
under the projective flow, and the same holds infinitesimally
on $\widehat M$.

The isometry group of the compact Riemannian manifold
$(\widehat M,h)$ is a compact Lie group.  The closure of the
descended one-parameter subgroup is therefore a
positive-dimensional compact connected abelian Lie group,
hence a torus.  Any circle subgroup of this torus gives
a nontrivial smooth $S^1$-action on $\widehat M$.
Yano's vanishing theorem \cite[p.~493]{Yano82} now gives
$\|\widehat M\|=0$.

Finally, multiplicativity of simplicial volume under finite
coverings yields
\[
 0=\|\widehat M\|=[\Gamma:\Gamma_0]\,\|M\|,
\]
so $\|M\|=0$.  This also covers nonorientable $M$, using the
orientation-cover convention adopted above.
\end{proof}


\begin{thebibliography}{99}
\bibitem[AMRT10]{AMRT10}
A.~Ash, D.~Mumford, M.~Rapoport, and Y.-S.~Tai,
\emph{Smooth Compactifications of Locally Symmetric Varieties},
second edition, Cambridge Mathematical Library,
Cambridge University Press, Cambridge, 2010.
\href{https://doi.org/10.1017/CBO9780511674693}
{doi:10.1017/CBO9780511674693}.

\bibitem[BDL18]{BallasDanciger}
S.~A.~Ballas, J.~Danciger, and G.-S.~Lee,
\emph{Convex projective structures on nonhyperbolic three-manifolds},
Geom. Topol. \textbf{22} (2018), no.~3, 1593--1646.
\href{https://doi.org/10.2140/gt.2018.22.1593}{doi:10.2140/gt.2018.22.1593}.


\bibitem[BP92]{BenedettiPetronio}
R.~Benedetti and C.~Petronio,
\emph{Lectures on Hyperbolic Geometry},
Universitext, Springer-Verlag, Berlin, 1992.
\href{https://doi.org/10.1007/978-3-642-58158-8}{doi:10.1007/978-3-642-58158-8}.

\bibitem[Ben08]{Benoist08Survey}
Y.~Benoist,
\emph{A survey on divisible convex sets},
Advanced Lectures in Mathematics \textbf{6} (2008), 1--18.


\bibitem[Buc07]{Bucher07}
M.~Bucher--Karlsson,
\emph{Simplicial volume of locally symmetric spaces covered by
$\mathrm{SL}_3\R/\mathrm{SO}(3)$},
Geom. Dedicata \textbf{125} (2007), 203--224.

\bibitem[BFR26]{BucherFrancavigliaRago26}
M.~Bucher, S.~Francaviglia, and G.~Rago,
\emph{Angle cocycles, Euler characteristic, and simplicial volume of
smooth manifolds},
\href{https://arxiv.org/abs/2610.32342}{arXiv:2610.32342v1}.  

\bibitem[Cal72]{Calabi72}
E.~Calabi,
\emph{Complete affine hyperspheres. I},
in \emph{Symposia Mathematica}, Vol.~X
(Convegno di Geometria Differenziale, INDAM, Rome, 1971),
Academic Press, London, 1972, pp.~19--38.


\bibitem[Fai21]{Faifman21}
D.~Faifman,
\emph{A Funk perspective on billiards, projective geometry and Mahler volume},
J. Diff. Geom. \textbf{127} (2024), no. 1, 161–212.
\href{https://doi.org/10.4310/jdg/1717356157}{doi:10.4310/jdg/1717356157}.

\bibitem[Fra18]{Fra18}
F. Franceschini,
\emph{A characterization of relatively hyperbolic groups via bounded cohomology},
Groups Geom. Dyn. \textbf{12} (2018), no.~3, 919--960.
\href{https://doi.org/10.4171/GGD/463}{doi:10.4171/GGD/463}.

\bibitem[GP21]{GoffiPediconi}
A.~Goffi and F.~Pediconi,
\emph{A note on the strong maximum principle for fully nonlinear
equations on Riemannian manifolds},
J. Geom. Anal. \textbf{31} (2021), 8641--8665.
\href{https://doi.org/10.1007/s12220-021-00607-2}{doi:10.1007/s12220-021-00607-2}.

\bibitem[Gro82]{Gromov82}
M.~Gromov,
\emph{Volume and bounded cohomology},
Publ. Math. Inst. Hautes \'Etudes Sci. \textbf{56} (1982), 5--99.
\href{https://www.numdam.org/item/PMIHES_1982__56__5_0/}{Numdam}.

\bibitem[HM81]{HM}
U.~Haagerup and H.~J.~Munkholm,
\emph{Simplices of maximal volume in hyperbolic $n$-space},
Acta Math. \textbf{147} (1981), 1--11.
\href{https://doi.org/10.1007/BF02392865}{doi:10.1007/BF02392865}.

\bibitem[Hil13]{Hildebrand13}
R.~Hildebrand,
\emph{On the $\infty$-norm of the cubic form of complete hyperbolic affine hyperspheres},
Results Math. \textbf{64} (2013), 113--119.
\href{https://doi.org/10.1007/s00025-012-0301-7}{doi:10.1007/s00025-012-0301-7}.

\bibitem[HLV08]{HuLiVrancken08}
Z.~Hu, H.~Li, and L.~Vrancken,
\emph{Characterizations of the Calabi product of hyperbolic affine hyperspheres},
Results Math. \textbf{52} (2008), 299--314.
\href{https://doi.org/10.1007/s00025-008-0312-6}{doi:10.1007/s00025-008-0312-6}.


\bibitem[Joh48]{John48}
F.~John,
\emph{Extremum problems with inequalities as subsidiary conditions},
in \emph{Studies and Essays Presented to R.~Courant on His 60th Birthday},
Interscience, New York, 1948, pp.~187--204.

\bibitem[KK17]{KlartagKolesnikov17}
B.~Klartag and A.~V.~Kolesnikov,
\emph{Remarks on curvature in the transportation metric},
Anal. Math. \textbf{43} (2017), 67--88.
\href{https://arxiv.org/abs/1604.04165}{arXiv:1604.04165}.

\bibitem[LS06]{LafontSchmidt06}
J.-F.~Lafont and B.~Schmidt,
\emph{Simplicial volume of closed locally symmetric spaces of non-compact type},
Acta Math. \textbf{197} (2006), 129--143.

\bibitem[LLS04]{LLS04}
A.-M. Li, H. Li, and U. Simon,
\emph{Centroaffine Bernstein problems},
Differential Geom. Appl. \textbf{20} (2004), no.~3, 331--356.
\href{https://doi.org/10.1016/j.difgeo.2003.12.004}{doi:10.1016/j.difgeo.2003.12.004}.

\bibitem[Lof08]{Loftin08}
J.~Loftin,
\emph{Survey on affine spheres},
Handbook of Geometric Analysis, no. 2 (2010).

\bibitem[Min01]{Mineyev01}
I.~Mineyev,
\emph{Straightening and bounded cohomology of hyperbolic groups},
Geom. Funct. Anal. \textbf{11} (2001), 807--839.

\bibitem[RSS26]{RioloSeppiSlavich26}
S.~Riolo, A.~Seppi, and L.~Slavich,
\emph{Convex projective manifolds, symmetric spaces and geometric
decompositions},
\href{https://arxiv.org/abs/2607.05662}{arXiv:2607.05662}.

\bibitem[RS57]{RS57}
C.~A.~Rogers and G.~C.~Shephard,
\emph{The difference body of a convex body},
Arch. Math. \textbf{8} (1957), 220--233.
doi:10.1007/BF01899997.

\bibitem[Yan82]{Yano82}
K.~Yano,
\emph{Gromov invariant and $S^1$-actions},
J. Fac. Sci. Univ. Tokyo Sect. IA Math.
\textbf{29} (1982), no.~3, 493--501.
\href{https://doi.org/10.15083/00039566}
{doi:10.15083/00039566}.


\end{thebibliography}
\end{document}